\documentclass[reqno,draft]{amsart}

\usepackage{amsmath}
\usepackage{amsthm}
\usepackage{amsfonts}
\usepackage{amssymb}
\usepackage{latexsym}
\usepackage{amscd}
\usepackage{stmaryrd}
\usepackage{color}

\theoremstyle{plain}
\newtheorem{thm}{Theorem}[section]
\newtheorem{cor}[thm]{Corollary}

\newtheorem{prop}[thm]{Proposition}
\newtheorem{rem}[thm]{Remark}

\numberwithin{equation}{section}

\newfont{\scyr}{wncyr10 scaled 550}
\def\shuffle{\,\mbox{\bf \scyr X}\,}

\def\Li{\operatorname{Li}}

\def\proof{\noindent {\bf Proof.\;}}

\allowdisplaybreaks

\begin{document}

\title{Shuffle product formulas of multiple zeta values}

\date{\today\thanks{The first author is supported by the National Natural Science Foundation of
China (Grant No. 11471245) and Shanghai Natural Science Foundation (grant no. 14ZR1443500). The authors are grateful to the referee for his/her useful remarks.} }

\author{Zhonghua Li \quad and \quad Chen Qin}

\address{Department of Mathematics, Tongji University, No. 1239 Siping Road,
Shanghai 200092, China}

\email{zhonghua\_li@tongji.edu.cn}

\address{Department of Mathematics, Tongji University, No. 1239 Siping Road,
Shanghai 200092, China}

\email{2014chen\_qin@tongji.edu.cn}



\begin{abstract}
Using the combinatorial description of shuffle product, we prove or reformulate several shuffle product formulas of multiple zeta values, including a general formula of the shuffle product of two multiple zeta values, some restricted shuffle product formulas of the product of two multiple zeta values, and a restricted shuffle product formula of the product of $n$ multiple zeta values.
\end{abstract}

\maketitle


\section{Introduction}\label{Sec:Intro}

For positive integers $n,k_1,k_2,\ldots,k_n$ with $k_1>1$, a multiple zeta value is the real number defined by
\begin{align}
\zeta(k_1,k_2,\ldots,k_n)=\sum\limits_{m_1>m_2>\cdots>m_n>0}\frac{1}{m_1^{k_1}m_2^{k_2}\cdots m_n^{k_n}}.
\label{Eq:MZV}
\end{align}
When $n=1$, we get the Riemann zeta values, which are the special values of the Riemann zeta function at positive integer arguments. There are many works on these real numbers. A recent theorem of F. Brown \cite{Brown} states that all periods of mixed Tate motives unramified over $\mathbb{Z}$ are $\mathbb{Q}\left[\frac{1}{2\pi i}\right]$-linear combinations of multiple zeta values, and the multiple zeta values indexed by $2$ and $3$ are linear generators of the $\mathbb{Q}$-vector space spanned by all multiple zeta values. In \cite{Minh2013}, Hoang Ngoc Minh showed that there exists also a family of algebraic generators, made up of the multiple zeta values indexed by irreducible Lyndon compositions.

Besides the infinite series representation \eqref{Eq:MZV}, N. Nielsen \cite{Nielsen04,Nielsen06} first noticed that a multiple zeta value can be also obtained via the iterated integral representation of the multiple polylogarithm, for $z$ tending to $1$, (see also \cite{Zagier}):
\begin{align}
\Li_{k_1,k_2,\ldots,k_n}(z)=\int\limits_{z>t_1>t_2>\cdots>t_k>0}\frac{dt_1}{f_1(t_1)}\frac{dt_2}{f_2(t_2)}\cdots\frac{dt_k}{f_k(t_k)},\quad \text{for\;} |z|<1,
\label{Eq:MZV-shuffle}
\end{align}
where $k=k_1+k_2+\cdots+k_n$ and
$$f_i(t)=\begin{cases}
1-t, & \text{if\;} i=k_1,k_1+k_2,\ldots,k_1+k_2+\cdots+k_{n-1},k,\\
t, & \text{otherwise}.
\end{cases}$$
Note that in \eqref{Eq:MZV-shuffle}, $k_1$ can be $1$. Using the iterated integral representation \eqref{Eq:MZV-shuffle}, we can express a product of two multiple zeta values as a sum of multiple zeta values. For example, we have
\begin{align*}
&\zeta(2)\zeta(2)=\int\limits_{1>t_1>t_2>0}\frac{dt_1}{t_1}\frac{dt_2}{1-t_2}\int\limits_{1>s_1>s_2>0}\frac{ds_1}{s_1}\frac{ds_2}{1-s_2}\\
=&\left(\,\int\limits_{1>t_1>t_2>s_1>s_2>0}+\int\limits_{1>t_1>s_1>t_2>s_2>0}+\int\limits_{1>t_1>s_1>s_2>t_2>0}\right.\\
&\left.+\int\limits_{1>s_1>t_1>t_2>s_2>0}+\int\limits_{1>s_1>t_1>s_2>t_2>0}+\int\limits_{1>s_1>s_2>t_1>t_2>0}\right)\\
&\quad\times \frac{dt_1}{t_1}\frac{dt_2}{1-t_2}\frac{ds_1}{s_1}\frac{ds_2}{1-s_2}\\
=&2\zeta(2,2)+4\zeta(3,1).
\end{align*}
Such products are called shuffle products. The shuffle product used here was first introduced by S. Eilenberg and S. Mac Lane in \cite{Eilenberg}, and the recursive formula and the denotation of the product described below are those of M. Fliess \cite{Fliess}.

To treat the shuffle products of multiple zeta values formally, we adopt the following algebraic setting (see \cite{Hoffman,Hoffman-Ohno} for example). Let $A=\{x,y\}$ be an alphabet with two non-commutative letters, and let $A^{\ast}$ be the set of all words on $A$ with the empty word $1_{A^{\ast}}$. We denote by $\mathfrak{h}$ the non-commutative $\mathbb{Q}$-polynomial algebra generated by the set $A$, and by $\mathfrak{h}^1$ and $\mathfrak{h}^0$ the subalgebras
$$\mathfrak{h}^1=\mathbb{Q}1_{A^{\ast}}+\mathfrak{h}y, \quad\mathfrak{h}^0=\mathbb{Q}1_{A^\ast}+x\mathfrak{h}y,$$
respectively. As rational vector spaces, $\mathfrak{h}$ is spanned by $A^{\ast}$, $\mathfrak{h}^1$ is spanned by $1_{A^{\ast}}$ and words ending with $y$ and $\mathfrak{h}^0$ is spanned by $1_{A^{\ast}}$ and words starting from $x$ and ending with $y$.
The shuffle product $\shuffle$ on $\mathfrak{h}$ is defined by $\mathbb{Q}$-bilinearity and the rules:
\begin{align*}
& 1_{A^{\ast}}\shuffle w=w\shuffle 1_{A^{\ast}}=w,\\
& aw_1\shuffle bw_2=a(w_1\shuffle bw_2)+b(aw_1\shuffle w_2),
\end{align*}
for all letters $a,b\in A$ and all words $w,w_1,w_2\in A^\ast$. Under the shuffle product, $\mathfrak{h}$ becomes a commutative $\mathbb{Q}$-algebra, and $\mathfrak{h}^1$ and $\mathfrak{h}^0$ are also subalgebras. As a commutative algebra, the shuffle algebra $\mathfrak{h}$ is free with a pure transcendence basis consisting of all Lyndon words \cite{Reutenauer}.

We define a $\mathbb{Q}$-linear map $\zeta:\mathfrak{h}^0\rightarrow \mathbb{R}$ by $\zeta(1_{A^{\ast}})=1$ and
\begin{align*}
\zeta(x^{k_1-1}yx^{k_2-1}y\cdots x^{k_n-1}y)=\zeta(k_1,\ldots,k_n),
\end{align*}
where $n,k_1,k_2,\ldots,k_n$ are positive integers with $k_1>1$. Then it is easy to know that the map $\zeta: (\mathfrak{h}^0,\shuffle)\rightarrow \mathbb{R}$ is an algebra homomorphism (see \cite{Hoffman-Ohno} for example). In other words, we have
$$\zeta(w_1\shuffle w_2)=\zeta(w_1)\zeta(w_2)$$
for any $w_1,w_2\in\mathfrak{h}^0$.

Note that for any $w\in\mathfrak{h}^1$, one can define the multiple polylogarithm $\Li_w(z)$ by $\mathbb{Q}$-linearities, $\Li_{1_{A^{\ast}}}(z)=1$ and
$$\Li_{x^{k_1-1}yx^{k_2-1}y\cdots x^{k_n-1}y}(z)=\Li_{k_1,k_2,\ldots,k_n}(z)$$
for positive integers $n,k_1,k_2,\ldots,k_n$. Then it is known that the map $\Li$ is an algebraic isomorphism from $(\mathfrak{h}^1,\shuffle)$ onto the algebra of multiple polylogarithms \cite{Minh-Petitot}.

Hence to get shuffle product formulas of multiple zeta values (or multiple polylogarithms), one can treat the shuffle products of $\mathfrak{h}^0$ first. For example, since
$$xy\shuffle xy=2xyxy+4x^2y^2,$$
applying the map $\zeta$, we get the shuffle product formula $\zeta(2)\zeta(2)=2\zeta(2,2)+4\zeta(3,1)$ proved above.

The shuffle product formula of two Riemann zeta values is
\begin{align}
\zeta(k)\zeta(l)=&\sum\limits_{i=1}^k\binom{k+l-i-1}{l-1}\zeta(k+l-i,i)\nonumber\\
&+\sum\limits_{i=1}^l\binom{k+l-i-1}{k-1}\zeta(k+l-i,i),
\label{Eq:Euler-Decom}
\end{align}
where $k,l\geqslant 2$. The formula \eqref{Eq:Euler-Decom} was first found by Euler \cite{Euler}, and is called Euler's decomposition formula. The corresponding Euler's decomposition formula in $\mathfrak{h}^1$ is
\begin{align}
x^{a}y\shuffle x^by=\sum\limits_{i=0}^a\binom{a+b-i}{b}x^{a+b-i}yx^iy+\sum\limits_{i=0}^b\binom{a+b-i}{a}x^{a+b-i}yx^iy,
\label{Eq:Shuffle-Euler}
\end{align}
where $a$ and $b$ are nonnegative integers. We remark that the shuffle product formula \eqref{Eq:Shuffle-Euler} is equivalent to
\begin{align*}
\Li_k(z)\Li_l(z)=&\sum\limits_{i=1}^k\binom{k+l-i-1}{l-1}\Li_{k+l-i,i}(z)\\
&+\sum\limits_{i=1}^l\binom{k+l-i-1}{k-1}\Li_{k+l-i,i}(z),
\end{align*}
where $k$ and $l$ are positive integers \cite{Minh-Petitot}.

Some generalizations of Euler's decomposition formula were found. In \cite[Theorem 2.1, Theorem 2.2]{Guo-Xie}, L. Guo and B. Xie gave an explicit shuffle product formula in a very general setting, and as applications, shuffle product formulas of $\zeta(k)\zeta(k_1,k_2)$ and $\zeta(k_1,k_2)\zeta(l_1,l_2)$ were given. By an analytic method, M. Eie and C.-S. Wei obtained a shuffle product formula of the product of two multiple zeta values of the form $\zeta(m,\{1\}^k)$ with one string of $1$'s: $\{1\}^k=\underbrace{1,\ldots,1}_{k\text{\, terms}}$ in \cite{Eie-Wei}. And this result was generalized to the product of $n$ multiple zeta values with one string of $1$'s in \cite[Main Theorem]{Eie-Liaw-Ong}, which also generalized the formula of the products of $n$ Riemann zeta values in \cite[Theorem 1.2]{Chung-Eie-Liaw-Ong}. By an algebraic method, P. Lie, L. Guo and B. Ma also gave a shuffle product formula of two multiple zeta values with  one string of $1$'s, and obtained a formula of the product of two multiple zeta values one of which with two strings of $1$'s in \cite[Theorem 1.1,Theorem 1.3]{Lei-Guo-Ma}.

The shuffle product has a combinatorial description. By the definition of shuffle product $\shuffle$ in $\mathfrak{h}$, we easily get
$$a_1\cdots a_n\shuffle a_{n+1}\cdots a_{n+m}=\sum\limits_{\sigma\in\mathfrak{S}_{n,m}}a_{\sigma(1)}a_{\sigma(2)}\cdots a_{\sigma(n+m)},$$
where $a_1, a_2,\ldots,a_{n+m}\in A$ are letters, and
$$\mathfrak{S}_{n,m}=\left\{\sigma\in\mathfrak{S}_{n+m}\left|\begin{array}{l}
\sigma^{-1}(1)<\sigma^{-1}(2)<\cdots<\sigma^{-1}(n),\\
\sigma^{-1}(n+1)<\sigma^{-1}(n+2)<\cdots<\sigma^{-1}(n+m)
\end{array}\right.\right\}.$$
In other words, the shuffle product of  words $a_1\cdots a_n$ and $a_{n+1}\cdots a_{n+m}$ are the sum of all permutations of $a_1,a_2,\ldots,a_{n+m}$, which simultaneously preserve the relative order of $a_1,a_2,\ldots,a_n$ and the relative order of $a_{n+1},a_{n+2},\ldots,a_{n+m}$. In this paper, we use this simple description to study shuffle products of words in $\mathfrak{h}y$, reformulate the formulas mentioned in the last paragraph, and find some new shuffle product formulas. We remark that in many cases, the method used here are very simple and natural. And the idea used here goes back to Hoang Ngoc Minh (see \cite{Minh96}-\cite{Minh-Petitot}).

We give the contexts of this paper. In Section \ref{Sec:Shu-MZV}, we give a general formula of the shuffle product of two words in $\mathfrak{h}y$ and provide some concrete examples. In Section \ref{Sec:Shu-MZV-1}, we give formulas of the shuffle products $x^ay^r\shuffle x^by^s$, $x^ay^r\shuffle x^{b_1}y^{s_1}x^{b_2}y^{s_2}$ and $x^{a_1}y^{r_1}x^{a_2}y^{r_2}\shuffle x^{b_1}y^{s_1}x^{b_2}y^{s_2}$. And we give a shuffle product formula of the products $x^{a_1}y^{r_1}\shuffle\cdots\shuffle x^{a_n}y^{r_n}$ in Section \ref{Sec:Shu-MZV-n}. There are two appendixes, in which we prove that the formulas found in \cite[Theorem 1.1,Theorem 1.3]{Lei-Guo-Ma} are essentially the same as \eqref{Eq:Shuffle-Res-1-1} and \eqref{Eq:Shuffle-Res-1-2}, respectively.


\section{Shuffle product formulas}\label{Sec:Shu-MZV}

\subsection{The general shuffle product formula}

As the most simple and typical example, we reprove Euler's decomposition formula \eqref{Eq:Shuffle-Euler} here.

Let $a$ and $b$ be two nonnegative integers. We want to compute the shuffle product $x^ay\shuffle x^by$. We may assume that
$$x^ay\shuffle x^by=\sum\limits_{\alpha_1+\alpha_2=a+b\atop \alpha_1,\alpha_2\geqslant 0}c_{\alpha_1,\alpha_2}x^{\alpha_1}yx^{\alpha_2}y.$$
For $\alpha_1,\alpha_2\geqslant 0$ with $\alpha_1+\alpha_2=a+b$, we have to determine the coefficients $c_{\alpha_1,\alpha_2}$. The ideal is that we shuffle the two $y$'s first, then consider where the $x$'s are from in $x^{\alpha_i}$. To distinguish these two $y$'s, we write the $y$ in $x^ay$ by $y_1$ and $y$ in $x^by$ by $y_2$. Then there are two cases when we shuffle $y_1$ and $y_2$:
\begin{description}
  \item[(i)] $y_1\quad y_2$;
  \item[(ii)] $y_2\quad y_1$.
\end{description}
For case (i), the number of $x$'s in $x^{\alpha_1}$ coming from $x^ay$ is $a$ and coming from $x^by$ is $\alpha_1-a$, respectively. And all $x$'s in $x^{\alpha_2}$ are coming from $x^by$. Hence the possibility is $\binom{\alpha_1}{a}$. Similarly, for case (ii), we get the possibility $\binom{\alpha_1}{b}$. Then we find the coefficient
$$c_{\alpha_1,\alpha_2}=\binom{\alpha_1}{a}+\binom{\alpha_1}{b}.$$
Finally we get the shuffle product formula
\begin{align}
x^ay\shuffle x^by=\sum\limits_{\alpha_1+\alpha_2=a+b\atop \alpha_1,\alpha_2\geqslant 0}\left[\binom{\alpha_1}{a}+\binom{\alpha_1}{b}\right]x^{\alpha_1}yx^{\alpha_2}y,
\label{Eq:Shuffle-Euler-New}
\end{align}
which is just Euler's decomposition formula \eqref{Eq:Shuffle-Euler}. In the summation of \eqref{Eq:Shuffle-Euler-New}, besides the condition $\alpha_1+\alpha_2=a+b$, the nonnegative integers $\alpha_1$ and $\alpha_2$ must satisfy the condition $\alpha_1\geqslant \min\{a,b\}$. While for integers $\alpha$ and $a$, the binomial coefficient $\binom{\alpha}{a}$ is zero if $a<0$ or $\alpha<a$. Therefore for simplicity, here and below we do not  write such type conditions.

Applying the same method to the general case, we get a general shuffle product formula, which is stated in the following theorem.

\begin{thm}\label{Thm:Shuffle-General}
Let $r,s$ be two positive integers and let $a_1,\ldots,a_r,b_1,\ldots,b_s$ be nonnegative integers. Then we have
\begin{align}
&x^{a_1}y\cdots x^{a_r}y\shuffle x^{b_1}y\cdots x^{b_s}y\nonumber\\
=&\sum\limits_{\alpha_1+\cdots+\alpha_{r+s}=\sum\limits_{i=1}^ra_i+\sum\limits_{j=1}^sb_j\atop \alpha_1,\ldots,\alpha_{r+s}\geqslant 0}c_{\alpha_1,\ldots,\alpha_{r+s}}x^{\alpha_1}yx^{\alpha_2}y\cdots x^{\alpha_{r+s}}y,
\label{Eq:Shuffle-General}
\end{align}
where the coefficients
\begin{align}
c_{\alpha_1,\ldots,\alpha_{r+s}}=&\sum\limits_{{l_1+\cdots+l_{p+1}=r\atop n_1+\cdots+n_{p}=s}\atop p\geqslant 1,l_i\geqslant 1,n_j\geqslant 1}\prod\limits_{i=1}^{L_p+s}\binom{\alpha_i}{\beta_i}\prod\limits_{j=L_p+s+2}^{r+s}\delta_{\alpha_j,a_{j-s}}\nonumber\\
&+\sum\limits_{{l_1+\cdots+l_{p}=r\atop n_1+\cdots+n_{p}=s}\atop p\geqslant 1,l_i\geqslant 1,n_j\geqslant 1}\prod\limits_{i=1}^{r+N_{p-1}}\binom{\alpha_i}{\beta_i}\prod\limits_{j=r+N_{p-1}+2}^{r+s}\delta_{\alpha_j,b_{j-r}}\nonumber\\
&+\sum\limits_{{l_1+\cdots+l_{p}=r\atop n_1+\cdots+n_{p+1}=s}\atop p\geqslant 1,l_i\geqslant 1,n_j\geqslant 1}\prod\limits_{i=1}^{r+N_p}\binom{\alpha_i}{\gamma_i}\prod\limits_{j=r+N_p+2}^{r+s}\delta_{\alpha_j,b_{j-r}}\nonumber\\
&+\sum\limits_{{l_1+\cdots+l_{p}=r\atop n_1+\cdots+n_{p}=s}\atop p\geqslant 1,l_i\geqslant 1,n_j\geqslant 1}\prod\limits_{i=1}^{L_{p-1}+s}\binom{\alpha_i}{\gamma_i}\prod\limits_{j=L_{p-1}+s+2}^{r+s}\delta_{\alpha_j,a_{j-s}}.
\label{Eq:Shuffle-General-Coe}
\end{align}
Here for positive integers $l_1,\ldots,l_p (,l_{p+1})$ and $n_1,\ldots,n_p (,n_{p+1})$ appearing in the summations above, we define
\begin{align}
\begin{cases}
\begin{array}{l}
\beta_{L_j+N_j+1}=\sum\limits_{i=1}^{L_j+1}a_i+\sum\limits_{i=1}^{N_j}b_i-\sum\limits_{i=1}^{L_j+N_j}\alpha_i, \\
\beta_{L_j+N_j+t}=a_{L_j+t}, \quad (2\leqslant t\leqslant l_{j+1})
\end{array} & \text{for\;} j=0,1,\ldots,p,\\
\begin{array}{l}
\beta_{L_{j+1}+N_j+1}=\sum\limits_{i=1}^{L_{j+1}}a_i+\sum\limits_{i=1}^{N_j+1}b_i-\sum\limits_{i=1}^{L_{j+1}+N_j}\alpha_i, \\
\beta_{L_{j+1}+N_j+t}=b_{N_j+t}, \quad (2\leqslant t\leqslant n_{j+1})
\end{array} & \text{for\;} j=0,1,\ldots,p-1,
\end{cases}
\label{Eq:Coe-Beta}
\end{align}
and
\begin{align}
\begin{cases}
\begin{array}{l}
\gamma_{L_j+N_j+1}=\sum\limits_{i=1}^{L_j}a_i+\sum\limits_{i=1}^{N_j+1}b_i-\sum\limits_{i=1}^{L_j+N_j}\alpha_i, \\
\gamma_{L_j+N_j+t}=b_{N_j+t}, \quad (2\leqslant t\leqslant n_{j+1})
\end{array} & \text{for\;} j=0,1,\ldots,p,\\
\begin{array}{l}
\gamma_{L_{j}+N_{j+1}+1}=\sum\limits_{i=1}^{L_{j}+1}a_i+\sum\limits_{i=1}^{N_{j+1}}b_i-\sum\limits_{i=1}^{L_{j}+N_{j+1}}\alpha_i, \\
\gamma_{L_{j}+N_{j+1}+t}=a_{L_j+t}, \quad (2\leqslant t\leqslant l_{j+1})
\end{array} & \text{for\;} j=0,1,\ldots,p-1,
\end{cases}
\label{Eq:Coe-Gamma}
\end{align}
with $L_j=l_1+\cdots+l_j$, $N_j=n_1+\cdots+n_j$ for $j\geqslant 0$ and $L_0=N_0=0$. And $\delta_{ij}$ is Kronecker's delta symbol defined as
$$\delta_{ij}=\begin{cases}
1, & \text{if\;} i=j,\\
0, & \text{otherwise}.
\end{cases}$$
\end{thm}

\proof We write the $y$ in $x^{a_1}y\cdots x^{a_r}y$ as $y_1$, and the $y$ in $x^{b_1}y\cdots x^{b_s}y$ as $y_2$. Then there are four cases when we shuffle $y_1$'s and $y_2$'s:
\begin{description}
  \item[(i)] $\underbrace{y_1 \cdots y_1}_{l_1} \underbrace{y_2\cdots y_2}_{n_1}\cdots \underbrace{y_1\cdots y_1}_{l_p}\underbrace{y_2\cdots y_2}_{n_p}\underbrace{y_1\cdots y_1}_{l_{p+1}}$,\\
      where $l_1+\cdots+l_{p+1}=r$, $n_1+\cdots+n_p=s$ with $p,l_i,n_j\geqslant 1$;
  \item[(ii)] $\underbrace{y_1\cdots y_1}_{l_1} \underbrace{y_2\cdots y_2}_{n_1}\cdots \underbrace{y_1\cdots y_1}_{l_p} \underbrace{y_2\cdots y_2}_{n_p}$,\\
  where $l_1+\cdots+l_p=r$, $n_1+\cdots+n_p=s$ with $p,l_i,n_j\geqslant 1$;
  \item[(iii)] $\underbrace{y_2\cdots y_2}_{n_1} \underbrace{y_1\cdots y_1}_{l_1} \cdots \underbrace{y_2\cdots y_2}_{n_p} \underbrace{y_1\cdots y_1}_{l_p}\underbrace{y_2\cdots y_2}_{n_{p+1}}$,\\
      where $l_1+\cdots+l_p=r$, $n_1+\cdots+n_{p+1}=s$ with $p,l_i,n_j\geqslant 1$;
  \item[(iv)] $\underbrace{y_2\cdots y_2}_{n_1} \underbrace{y_1\cdots y_1}_{l_1} \cdots \underbrace{y_2\cdots y_2}_{n_p} \underbrace{y_1\cdots y_1}_{l_p}$,\\
      where $l_1+\cdots+l_p=r$, $n_1+\cdots+n_{p}=s$ with $p,l_i,n_j\geqslant 1$.
\end{description}

We compute the contribution of each case to the coefficient $c_{\alpha_1,\ldots,\alpha_{r+s}}$. Then $c_{\alpha_1,\ldots,\alpha_{r+s}}$ is just the sum of all these contributions. For case (i), we define
$$\beta_{L_j+N_j+t}=\text{the number of\;} x\text{'s coming from\;} x^{a_{L_j+t}}y \text{\;in\;} x^{\alpha_{L_j+N_j+t}}$$
for $0\leqslant j\leqslant p$ and $1\leqslant t\leqslant l_{j+1}$, and
$$\beta_{L_{j+1}+N_j+t}=\text{the number of\;} x\text{'s coming from\;} x^{b_{N_j+t}}y \text{\;in\;} x^{\alpha_{L_{j+1}+N_j+t}}$$
for $0\leqslant j\leqslant p-1$ and $1\leqslant t\leqslant n_{j+1}$. Note that the $x$'s in $x^{\alpha_{L_j+N_j+t}}$ are either from $x^{a_{L_j+t}}y$ or from $x^{b_{N_{j}+1}}y$, and these $x$'s can be in any order. Similarly, the $x$'s in $x^{\alpha_{L_{j+1}+N_j+t}}$ are either from $x^{b_{N_j+t}}y$ or from $x^{a_{L_{j+2}+1}}y$, and these $x$'s can be in any order too. Hence the contribution of case (i) to the coefficient $c_{\alpha_1,\ldots,\alpha_{r+s}}$ is
$$\sum\limits_{{l_1+\cdots+l_{p+1}=r\atop n_1+\cdots+n_{p}=s}\atop p\geqslant 1,l_i\geqslant 1,n_j\geqslant 1}\prod\limits_{i=1}^{r+s}\binom{\alpha_i}{\beta_i}.$$

We prove that the $\beta_i$'s satisfy the formulas given in \eqref{Eq:Coe-Beta} by induction on $j$.
When $j=0$, it is obvious that
\begin{align*}
&\beta_1=a_1,\quad \beta_2=a_2,\quad \ldots, \quad \beta_{L_1}=a_{L_1},\\
&\beta_{L_1+1}=b_1-\sum\limits_{i=1}^{L_1}(\alpha_i-\beta_i)=\sum\limits_{i=1}^{L_1}a_i+b_1-\sum\limits_{i=1}^{L_1}\alpha_i,\\
&\beta_{L_1+2}=b_2,\quad \ldots,\quad \beta_{L_1+n_1}=b_{n_1}.
\end{align*}
Now assume that $j>0$, and the formulas hold for $j-1$. Then we have
\begin{align*}
&\beta_{L_j+N_j+1}=a_{L_j+1}-\sum\limits_{i=L_j+N_{j-1}+1}^{L_j+N_j}(\alpha_i-\beta_i)\\
=&a_{L_j+1}+\sum\limits_{i=1}^{L_j}a_i+\sum\limits_{i=1}^{N_{j-1}+1}b_i-\sum\limits_{i=1}^{L_j+N_{j-1}}\alpha_i+\sum\limits_{t=2}^{n_j}b_{N_{j-1}+t}
-\sum\limits_{i=L_j+N_{j-1}+1}^{L_j+N_j}\alpha_i\\
=&\sum\limits_{i=1}^{L_j+1}a_i+\sum\limits_{i=1}^{N_{j}}b_i-\sum\limits_{i=1}^{L_j+N_{j}}\alpha_i.
\end{align*}
For $2\leqslant t\leqslant l_{j+1}$, it is obvious that $\beta_{L_j+N_j+t}=a_{L_j+t}$. Next we have
\begin{align*}
&\beta_{L_{j+1}+N_j+1}=b_{N_j+1}-\sum\limits_{i=L_j+N_j+1}^{L_{j+1}+N_j}(\alpha_i-\beta_i)\\
=&b_{N_j+1}+\sum\limits_{i=1}^{L_j+1}a_i+\sum\limits_{i=1}^{N_{j}}b_i-\sum\limits_{i=1}^{L_j+N_{j}}\alpha_i+\sum\limits_{t=2}^{l_{j+1}}a_{L_j+t}
-\sum\limits_{i=L_j+N_j+1}^{L_{j+1}+N_j}\alpha_i\\
=&\sum\limits_{i=1}^{L_{j+1}}a_i+\sum\limits_{i=1}^{N_j+1}b_i-\sum\limits_{i=1}^{L_{j+1}+N_j}\alpha_i.
\end{align*}
Finally, it is obvious that $\beta_{L_{j+1}+N_j+t}=b_{N_j+t}$ for $2\leqslant t\leqslant n_{j+1}$. Then we have proved that the $\beta_i$'s defined above are given by \eqref{Eq:Coe-Beta}.

It is easy to see that for $j$ with the conditions $L_p+N_p+2\leqslant j\leqslant r+s$, it must hold $\alpha_j=a_{j-N_p}$. We also know that $\alpha_{L_p+N_p+1}=\beta_{L_p+N_p+1}$. Hence the contribution of case (i) to the coefficient $c_{\alpha_1,\ldots,\alpha_{r+s}}$ is
$$\sum\limits_{{l_1+\cdots+l_{p+1}=r\atop n_1+\cdots+n_{p}=s}\atop p\geqslant 1,l_i\geqslant 1,n_j\geqslant 1}\prod\limits_{i=1}^{L_p+s}\binom{\alpha_i}{\beta_i}\prod\limits_{j=L_p+s+2}^{r+s}\delta_{\alpha_j,a_{j-s}}$$
with $\beta_i$'s given by \eqref{Eq:Coe-Beta}.

Similar to case (i), we find that the contribution of case (ii) to the coefficient $c_{\alpha_1,\ldots,\alpha_{r+s}}$ is
$$\sum\limits_{{l_1+\cdots+l_{p}=r\atop n_1+\cdots+n_{p}=s}\atop p\geqslant 1,l_i\geqslant 1,n_j\geqslant 1}\prod\limits_{i=1}^{r+N_{p-1}}\binom{\alpha_i}{\beta_i}\prod\limits_{j=r+N_{p-1}+2}^{r+s}\delta_{\alpha_j,b_{j-r}},$$
where the $\beta_i$'s are also given by \eqref{Eq:Coe-Beta}. Finally by symmetry, we get the contribution of case (iii) from case (i), and the contribution of case (iv) from case (ii). Thus the theorem is proved.
\qed

Applying the algebra homomorphism $\zeta:(\mathfrak{h}^0,\shuffle)\rightarrow \mathbb{R}$, we get the shuffle product formulas of multiple zeta values.

\begin{cor}
Let $r,s$ be two positive integers and let $a_1,\ldots,a_r,b_1,\ldots,b_s$ be nonnegative integers with $a_1,b_1\geqslant 1$. Then we have
\begin{align}
&\zeta(a_1+1,\ldots, a_r+1)\zeta(b_1+1,\ldots,b_s+1)\nonumber\\
=&\sum\limits_{\alpha_1+\cdots+\alpha_{r+s}=\sum\limits_{i=1}^ra_i+\sum\limits_{j=1}^sb_j\atop \alpha_1\geqslant 1,\alpha_2,\ldots,\alpha_{r+s}\geqslant 0}c_{\alpha_1,\ldots,\alpha_{r+s}}\zeta(\alpha_1+1,\alpha_2+1,\ldots, \alpha_{r+s}+1),
\label{Eq:MZV-Shuffle-General}
\end{align}
where the coefficients $c_{\alpha_1,\ldots,\alpha_{r+s}}$ are given by \eqref{Eq:Shuffle-General-Coe}.
\end{cor}

\begin{rem}
The shuffle product formula \eqref{Eq:MZV-Shuffle-General} is essentially the same as that of \cite[Corollary 2.5]{Guo-Xie}. The proof we supply here seems much more simple. In fact, one can prove \cite[Equaiton (22) in Theorem 2.1]{Guo-Xie} by the method provided in this paper.
\end{rem}

\begin{rem}
 Applying the isomorphism $\Li$ to the shuffle product formula \eqref{Eq:MZV-Shuffle-General} in Theorem \ref{Thm:Shuffle-General}, one can obtain an equivalent formula which holds for multiple polylogarithms.
\end{rem}

\begin{rem}
The quantities $c_{\alpha_1,\ldots,\alpha_{r+s}}$ appearing in \eqref{Eq:MZV-Shuffle-General} are the structure constants of the shuffle algebra $\mathfrak{h}$, obtained when one identifies the coefficients in the rational expressions \cite{Berstel} $(x+y)^{\ast}\shuffle(x+y)^\ast=(2x+2y)^{\ast}$.
\end{rem}

\begin{rem}
There is another commutative product among multiple zeta values, which is called stuffle product \cite{Hoffman}. If one can get the formula for the stuffle product
$x^{a_1}y\cdots x^{a_r}y\ast x^{b_1}y\cdots x^{b_s}y$ (we consider this product in \cite{Li-Qin}), then comparing with \eqref{Eq:Shuffle-General}, one can get some double shuffle relations of multiple zeta values.
While it seems that the obtained double shuffle relations are not terse and elegant. For example, the formula deduced from the difference $x^ay\shuffle x^{b_1}y\cdots x^{b_s}y-x^ay\ast x^{b_1}y\cdots x^{b_s}y$ is not so terse (see Proposition \ref{Prop:Shu-1-s} below). On the other hand, if one considers the difference of sums
$$\sum\limits_{a+b_1+\cdots+b_s=k}x^ay\shuffle x^{b_1}y\cdots x^{b_s}y,\quad \sum\limits_{a+b_1+\cdots+b_s=k}x^ay\ast x^{b_1}y\cdots x^{b_s}y,$$
one obtains an elegant formula, which is called weighted sum formula (see \cite[Theorems 2.2-2.6]{Guo-Xie-09}).
\end{rem}

\subsection{Some concrete examples}

We give some applications of Theorem \ref{Thm:Shuffle-General}. First, we consider the case $r=1$ and get the following result.

\begin{prop}\label{Prop:Shu-1-s}
Let $s$ be a positive integer and let $a,b_1,\ldots,b_s$ be nonnegative integers. Then we have
\begin{align}
&x^ay\shuffle x^{b_1}y\cdots x^{b_s}y=\sum\limits_{\alpha_1+\cdots+\alpha_{s+1}=a+b_1+\cdots+b_s\atop \alpha_1,\ldots,\alpha_{s+1}\geqslant 0}\left\{\binom{\alpha_1}{a}\prod\limits_{j=3}^{s+1}\delta_{\alpha_j,b_{j-1}}+\prod\limits_{i=1}^s\binom{\alpha_i}{b_i}\right.\nonumber\\
&\quad +\left.\sum\limits_{k=1}^{s-1}\prod\limits_{i=1}^{k}\binom{\alpha_i}{b_i}\binom{\alpha_{k+1}}{b_{k+1}-\alpha_{k+2}}\prod\limits_{j=k+3}^{s+1}
\delta_{\alpha_j,b_{j-1}}\right\}x^{\alpha_1}yx^{\alpha_2}y\cdots x^{\alpha_{s+1}}y.
\label{Eq:Shuffle-1-s}
\end{align}
If further $a,b_1\geqslant 1$, then we have
\begin{align}
&\zeta(a+1)\zeta(b_1+1,\ldots, b_s+1)=\sum\limits_{\alpha_1+\cdots+\alpha_{s+1}=a+b_1+\cdots+b_s\atop \alpha_1\geqslant 1,\alpha_2,\ldots,\alpha_{s+1}\geqslant 0}\left\{\binom{\alpha_1}{a}\prod\limits_{j=3}^{s+1}\delta_{\alpha_j,b_{j-1}}\right.\nonumber\\
&\quad \left. +\prod\limits_{i=1}^s\binom{\alpha_i}{b_i}+\sum\limits_{k=1}^{s-1}\prod\limits_{i=1}^{k}\binom{\alpha_i}{b_i}\binom{\alpha_{k+1}}{b_{k+1}-\alpha_{k+2}}\prod\limits_{j=k+3}^{s+1}
\delta_{\alpha_j,b_{j-1}}\right\}\nonumber\\
&\qquad \times \zeta(\alpha_1+1,\alpha_2+1,\ldots, \alpha_{s+1}+1).
\label{Eq:MZV-Shuffle-1-s}
\end{align}
\end{prop}

\proof Applying the map $\zeta$, we get \eqref{Eq:MZV-Shuffle-1-s} from \eqref{Eq:Shuffle-1-s}. Hence we only need to prove \eqref{Eq:Shuffle-1-s}. By Theorem \ref{Thm:Shuffle-General}, we have
$$x^ay\shuffle x^{b_1}y\cdots x^{b_s}y=\sum\limits_{\alpha_1+\cdots+\alpha_{s+1}=a+b_1+\cdots+b_s\atop \alpha_1,\ldots,\alpha_{s+1}\geqslant 0}c_{\alpha_1,\ldots,\alpha_{s+1}}x^{\alpha_1}yx^{\alpha_2}y\cdots x^{\alpha_{s+1}}y,$$
where
\begin{align*}
c_{\alpha_1,\ldots,\alpha_{s+1}}=&\sum\limits_{l_1=r,n_1=s}\binom{\alpha_1}{\beta_1}\prod\limits_{j=3}^{s+1}\delta_{\alpha_j,b_{j-1}}
+\sum\limits_{l_1=r,n_1+n_2=s\atop n_1,n_2\geqslant 1}\prod\limits_{i=1}^{n_1+1}\binom{\alpha_i}{\gamma_i}\prod\limits_{j=n_1+3}^{s+1}\delta_{\alpha_j,b_{j-1}}\\
&+\sum\limits_{l_1=r,n_1=s}\prod\limits_{i=1}^s\binom{\alpha_i}{\gamma_i}.
\end{align*}
Now for $l_1=r,n_1=s$, we have $\beta_1=a$ and
$$\gamma_i=b_i,\quad (i=1,\ldots,n_1).$$
And for $l_1=r, n_1+n_2=s$, we have
$$\gamma_i=b_i,\quad (i=1,\ldots,n_1),$$
and
$$\gamma_{n_1+1}=a+\sum\limits_{i=1}^{n_1}b_i-\sum\limits_{i=1}^{n_1}\alpha_i.$$
Then we get the formula \eqref{Eq:Shuffle-1-s}.
\qed

Let $s=1$ and $b_1=b$ in \eqref{Eq:Shuffle-1-s}, we get Euler's decomposition formula \eqref{Eq:Shuffle-Euler-New}. Let $s=2$, we have
\begin{align}
x^ay\shuffle x^{b_1}yx^{b_2}y=\sum\limits_{\alpha_1+\alpha_2+\alpha_3=a+b_1+b_2\atop \alpha_1,\alpha_2,\alpha_3\geqslant 0}\left\{\binom{\alpha_1}{a}\delta_{\alpha_3,b_2}+\binom{\alpha_1}{b_1}\binom{\alpha_2}{b_2}\right.\nonumber\\
\left.+\binom{\alpha_1}{b_1}\binom{\alpha_2}{b_2-\alpha_3}\right\}x^{\alpha_1}yx^{\alpha_2}yx^{\alpha_3}y.
\label{Eq:Shuffle-1-2}
\end{align}
Under the conditions $a,b_1\geqslant 1$, after applying the map $\zeta$ to \eqref{Eq:Shuffle-1-2}, we get \cite[Equation (3)]{Guo-Xie}.
Similarly, let $s=3$, we get
\begin{align}
&x^ay\shuffle x^{b_1}yx^{b_2}yx^{b_3}y=\sum\limits_{{\alpha_1+\alpha_2+\alpha_3+\alpha_4\atop =a+b_1+b_2+b_3}\atop \alpha_1,\ldots,\alpha_4\geqslant 0}\left\{\binom{\alpha_1}{a}\delta_{\alpha_3,b_2}\delta_{\alpha_4,b_3}
+\binom{\alpha_1}{b_1}\binom{\alpha_2}{b_2-\alpha_3}\delta_{\alpha_4,b_3}\right.\nonumber\\
&\qquad \left.+\binom{\alpha_1}{b_1}\binom{\alpha_2}{b_2}
\left[\binom{\alpha_3}{b_3}+\binom{\alpha_3}{b_3-\alpha_4}\right]\right\}x^{\alpha_1}yx^{\alpha_2}yx^{\alpha_3}yx^{\alpha_4}y.
\label{Eq:Shuffle-1-3}
\end{align}

Next we compute the formula of $r=s=2$. By Theorem \ref{Thm:Shuffle-General}, we have
$$x^{a_1}yx^{a_2}y\shuffle x^{b_1}yx^{b_2}y=\sum\limits_{{\alpha_1+\alpha_2+\alpha_3+\alpha_4\atop=a_1+a_2+b_1+b_2}\atop \alpha_1,\ldots,\alpha_4\geqslant 0}(c_1+c_2+c_1'+c_2') x^{\alpha_1}yx^{\alpha_2}yx^{\alpha_3}yx^{\alpha_4}y,$$
where
\begin{align*}
c_1=&\sum\limits_{{l_1=l_{2}=1,n_1=2}}\prod\limits_{i=1}^{3}\binom{\alpha_i}{\beta_i}=\binom{\alpha_1}{a_1}\binom{\alpha_2}{a_1+b_1-\alpha_1}\binom{\alpha_3}{b_2},\\
c_2=&\sum\limits_{{l_1=2,n_1=2}}\prod\limits_{i=1}^{2}\binom{\alpha_i}{\beta_i}\delta_{\alpha_4,b_{2}}+\sum\limits_{{l_1=l_{2}=1,n_1=n_{2}=1}}
\prod\limits_{i=1}^{3}\binom{\alpha_i}{\beta_i}\\
=&\binom{\alpha_1}{a_1}\binom{\alpha_2}{a_2}\delta_{\alpha_4,b_2}+\binom{\alpha_1}{a_1}\binom{\alpha_2}{a_1+b_1-\alpha_1}\binom{\alpha_3}{a_1+a_2+b_1-\alpha_1-\alpha_2}\\
=&\binom{\alpha_1}{a_1}\binom{\alpha_2}{a_2}\delta_{\alpha_4,b_2}+\binom{\alpha_1}{a_1}\binom{\alpha_2}{a_1+b_1-\alpha_1}\binom{\alpha_3}{b_2-\alpha_4},
\end{align*}
and $c_i'$ is obtained from $c_i$ by interchanging $a_1$ with $b_1$ and $a_2$ with $b_2$ simultaneously for $i=1,2$. Then we get the shuffle product formula
\begin{align}
&x^{a_1}yx^{a_2}y\shuffle x^{b_1}yx^{b_2}y=\sum\limits_{{\alpha_1+\alpha_2+\alpha_3+\alpha_4\atop=a_1+a_2+b_1+b_2}\atop \alpha_1,\ldots,\alpha_4\geqslant 0}\left\{\binom{\alpha_1}{a_1}\binom{\alpha_2}{a_2}\delta_{\alpha_4,b_2}+\binom{\alpha_1}{b_1}\binom{\alpha_2}{b_2}\delta_{\alpha_4,a_2}\right.\nonumber\\
&\quad +\binom{\alpha_1}{a_1}\binom{\alpha_2}{a_1+b_1-\alpha_1}\left[\binom{\alpha_3}{b_2}+\binom{\alpha_3}{b_2-\alpha_4}\right]\nonumber\\
&\quad\left.+\binom{\alpha_1}{b_1}\binom{\alpha_2}{a_1+b_1-\alpha_1}\left[\binom{\alpha_3}{a_2}+\binom{\alpha_3}{a_2-\alpha_4}\right]\right\}x^{\alpha_1}yx^{\alpha_2}yx^{\alpha_3}yx^{\alpha_4}y.
\label{Eq:Shuffle-2-2}
\end{align}
Under the conditions $a_1,b_1\geqslant 1$, after applying the map $\zeta$ to \eqref{Eq:Shuffle-2-2}, we get \cite[Equation (4)]{Guo-Xie}.

Finally, it is not difficult but needs patience to find the formulas of $r=2,s=3$ and $r=s=3$. The shuffle product formulas are
\begin{align}
&x^{a_1}yx^{a_2}y\shuffle x^{b_1}yx^{b_2}yx^{b_3}y=\sum\limits_{{\alpha_1+\cdots+\alpha_5\atop =a_1+a_2+b_1+b_2+b_3}\atop \alpha_1,\ldots,\alpha_5\geqslant 0}\left\{\binom{\alpha_1}{a_1}\binom{\alpha_2}{a_2}\delta_{\alpha_4,b_2}\delta_{\alpha_5,b_3}\right.\nonumber\\
&\quad+\binom{\alpha_1}{a_1}\binom{\alpha_2}{a_1+b_1-\alpha_1}\binom{\alpha_3}{b_2-a_4}\delta_{\alpha_5,b_3}\nonumber\\
&\quad+\binom{\alpha_1}{a_1}\binom{\alpha_2}{a_1+b_1-\alpha_1}\binom{\alpha_3}{b_2}\left[\binom{\alpha_4}{b_3}+\binom{\alpha_4}{b_3-\alpha_5}\right]\nonumber\\
&\quad+\binom{\alpha_1}{b_1}\binom{\alpha_2}{b_2}\binom{\alpha_3}{b_3}\delta_{\alpha_5,a_2}
+\binom{\alpha_1}{b_1}\binom{\alpha_2}{a_1+b_1-\alpha_1}\binom{\alpha_3}{a_2}\delta_{\alpha_5,b_3}\nonumber\\
&\quad+\binom{\alpha_1}{b_1}\binom{\alpha_2}{b_2}\binom{\alpha_3}{a_2+b_3-\alpha_4-\alpha_5}\left[\binom{\alpha_4}{a_2}+\binom{\alpha_4}{a_2-\alpha_5}\right]
\nonumber\\
&\quad\left.+\binom{\alpha_1}{b_1}\binom{\alpha_2}{a_1+b_1-\alpha_1}\binom{\alpha_3}{a_2+b_3-\alpha_4-\alpha_5}\left[\binom{\alpha_4}{b_3}
+\binom{\alpha_4}{b_3-\alpha_5}\right]\right\}\nonumber\\
&\qquad \times x^{\alpha_1}yx^{\alpha_2}yx^{\alpha_3}yx^{\alpha_4}yx^{\alpha_5}y,
\label{Eq:Shuffle-2-3}
\end{align}
and
\begin{align}
&x^{a_1}yx^{a_2}yx^{a_3}y\shuffle x^{b_1}yx^{b_2}yx^{b_3}y=\sum\limits_{{\alpha_1+\cdots+\alpha_6\atop =a_1+a_2+a_3+b_1+b_2+b_3}\atop \alpha_1,\ldots,\alpha_6\geqslant 0}[c(a_1,a_2;b_1,b_2)\nonumber\\
&\qquad \qquad +c(b_1,b_2;a_1,a_2)]x^{\alpha_1}yx^{\alpha_2}yx^{\alpha_3}yx^{\alpha_4}yx^{\alpha_5}yx^{\alpha_6}y,
\label{Eq:Shuffle-3-3}
\end{align}
with
\begin{align*}
&c(a_1,a_2;b_1,b_2)\\
=&\binom{\alpha_1}{a_1}\binom{\alpha_2}{a_2}\binom{\alpha_3}{a_3}\delta_{\alpha_5,b_2}\delta_{\alpha_6,b_3}+\binom{\alpha_1}{a_1}
\binom{\alpha_2}{a_1+b_1-\alpha_1}\binom{\alpha_3}{b_2}\binom{\alpha_4}{b_3}\delta_{\alpha_6,a_3}\\
&+\binom{\alpha_1}{a_1}\binom{\alpha_2}{a_1+b_1-\alpha_1}\binom{\alpha_3}{a_1+a_2+b_1-\alpha_1-\alpha_2}\binom{\alpha_4}{a_3}\delta_{\alpha_6,b_3}\\
&+\binom{\alpha_1}{a_1}\binom{\alpha_2}{a_2}\binom{\alpha_3}{a_1+a_2+b_1-\alpha_1-\alpha_2}\binom{\alpha_4}{b_2-\alpha_5}\delta_{\alpha_6,b_3}\\
&+\binom{\alpha_1}{a_1}\binom{\alpha_2}{a_2}\binom{\alpha_3}{a_1+a_2+b_1-\alpha_1-\alpha_2}\binom{\alpha_4}{b_2}
\left[\binom{\alpha_5}{b_3}+\binom{\alpha_5}{b_3-\alpha_6}\right]\\
&+\binom{\alpha_1}{a_1}\binom{\alpha_2}{a_1+b_1-\alpha_1}\binom{\alpha_3}{b_2}\binom{\alpha_4}{a_3+b_3-\alpha_5-\alpha_6}
\left[\binom{\alpha_5}{a_3}+\binom{\alpha_5}{a_3-\alpha_6}\right]\\
&+\binom{\alpha_1}{a_1}\binom{\alpha_2}{a_1+b_1-\alpha_1}\binom{\alpha_3}{a_1+a_2+b_1-\alpha_1-\alpha_2}\binom{\alpha_4}{a_3+b_3-\alpha_5-\alpha_6}\\
&\quad \times \left[\binom{\alpha_5}{b_3}+\binom{\alpha_5}{b_3-\alpha_6}\right].
\end{align*}
We omit the proofs of \eqref{Eq:Shuffle-2-3} and \eqref{Eq:Shuffle-3-3}. Of course, applying the map $\zeta$, we get shuffle product formulas of multiple zeta values from \eqref{Eq:Shuffle-2-3} and \eqref{Eq:Shuffle-3-3}.


\section{Restricted shuffle product formulas}\label{Sec:Shu-MZV-1}

 In this section, we study shuffle product formulas of multiple zeta values with strings of $1$'s. Such formulas are called restricted shuffle product formulas.

\subsection{The formula of $x^ay^r\shuffle x^by^s$}

In this subsection, we consider the shuffle product of multiple zeta values of the form $\zeta(m,\{1\}^n)$, which had been studied in \cite{Eie-Wei,Lei-Guo-Ma}.

Equivalently, we have to compute $x^ay^r\shuffle x^by^s$ for any nonnegative integers $a,b$ and any positive integers $r,s$. For that purpose, we can use Theorem \ref{Thm:Shuffle-General}. We have
$$x^ay^r\shuffle x^by^s=\sum\limits_{\alpha_1+\cdots+\alpha_{r+s}=a+b\atop \alpha_1,\ldots,\alpha_{r+s}\geqslant 0}[c(a,r;b,s)+c(b,s;a,r)]x^{\alpha_1}y\cdots x^{\alpha_{r+s}}y,$$
where $c(a,r;b,s)=\Sigma_1+\Sigma_2$ with
\begin{align*}
\Sigma_1=&\sum\limits_{{l_1+\cdots+l_{p+1}=r\atop n_1+\cdots+n_{p}=s}\atop p\geqslant 1,l_i\geqslant 1,n_j\geqslant 1}\prod\limits_{i=1}^{L_p+s}\binom{\alpha_i}{\beta_i}\prod\limits_{j=L_p+s+2}^{r+s}\delta_{\alpha_j,a_{j-s}},\\
\Sigma_2=&\sum\limits_{{l_1+\cdots+l_{p}=r\atop n_1+\cdots+n_{p}=s}\atop p\geqslant 1,l_i\geqslant 1,n_j\geqslant 1}\prod\limits_{i=1}^{r+N_{p-1}}\binom{\alpha_i}{\beta_i}\prod\limits_{j=r+N_{p-1}+2}^{r+s}\delta_{\alpha_j,b_{j-r}}.
\end{align*}
For $\Sigma_1$, we have
$$\begin{cases}
\beta_1=a, & \\
\beta_{L_j+N_j+1}=a+b-\sum\limits_{i=1}^{L_j+N_j}\alpha_i, & (j=1,\ldots,p),\\
\beta_{L_j+N_j+t}=0, & (2\leqslant t\leqslant l_{j+1},j=0,1,\ldots,p),\\
\beta_{L_{j+1}+N_j+1}=a+b-\sum\limits_{i=1}^{L_{j+1}+N_j}\alpha_i, & (j=0,1,\ldots,p-1),\\
\beta_{L_{j+1}+N_j+t}=0, & (2\leqslant t\leqslant n_{j+1},j=0,1,\ldots,p-1).
\end{cases}$$
Since $\alpha_{l_1+1}\geqslant \beta_{l_1+1}$, we get $\sum\limits_{i=1}^{l_1+1}\alpha_i\geqslant a+b$. While $a+b=\sum\limits_{i=1}^{r+s}\alpha_i$, hence it must hold $\alpha_{l_1+2}=\cdots=\alpha_{r+s}=0$ and $\sum\limits_{i=1}^{l_1+1}\alpha_i=a+b$. Then we get
\begin{align*}
\Sigma_1=&\sum\limits_{{l_1+\cdots+l_{p+1}=r\atop n_1+\cdots+n_{p}=s}\atop p\geqslant 1,l_i\geqslant 1,n_j\geqslant 1}\binom{\alpha_1}{a}\prod\limits_{j=l_1+2}^{r+s}\delta_{\alpha_j,0}\\
=&\sum\limits_{l=1}^{r-1}\left(\sum\limits_{{l_2+\cdots+l_{p+1}=r-l\atop n_1+\cdots+n_{p}=s}\atop p\geqslant 1,l_i\geqslant 1,n_j\geqslant 1}1\right)\binom{\alpha_1}{a}\prod\limits_{j=l+2}^{r+s}\delta_{\alpha_j,0}\\
=&\sum\limits_{l=1}^{r-1}\left(\sum\limits_{p\geqslant 1}\binom{r-l-1}{p-1}\binom{s-1}{p-1}\right)\binom{\alpha_1}{a}\prod\limits_{j=l+2}^{r+s}\delta_{\alpha_j,0}.
\end{align*}
Using the well-known combinatorial identity
\begin{align}
\sum\limits_{i=0}^n\binom{k}{i}\binom{l}{n-i}=\binom{k+l}{n},
\label{Eq:Com-Id}
\end{align}
we get
$$\Sigma_1=\sum\limits_{l=1}^{r-1}\binom{r+s-l-2}{r-l-1}\binom{\alpha_1}{a}\prod\limits_{j=l+2}^{r+s}\delta_{\alpha_j,0}.$$
Similarly, we have
\begin{align*}
\Sigma_2=&\sum\limits_{{l_1+\cdots+l_{p}=r\atop n_1+\cdots+n_{p}=s}\atop p\geqslant 1,l_i\geqslant 1,n_j\geqslant 1}\binom{\alpha_1}{a}\prod\limits_{j=l_1+2}^{r+s}\delta_{\alpha_j,0}\\
=&\sum\limits_{l=1}^{r-1}\left(\sum\limits_{{l_2+\cdots+l_{p}=r-l\atop n_1+\cdots+n_{p}=s}\atop p\geqslant 2,l_i\geqslant 1,n_j\geqslant 1}1\right)\binom{\alpha_1}{a}\prod\limits_{j=l+2}^{r+s}\delta_{\alpha_j,0}+\binom{\alpha_1}{a}\prod\limits_{j=r+2}^{r+s}\delta_{\alpha_j,0}\\
=&\sum\limits_{l=1}^{r-1}\left(\sum\limits_{p\geqslant 2}\binom{r-l-1}{p-2}\binom{s-1}{p-1}\right)\binom{\alpha_1}{a}\prod\limits_{j=l+2}^{r+s}\delta_{\alpha_j,0}+\binom{\alpha_1}{a}\prod\limits_{j=r+2}^{r+s}\delta_{\alpha_j,0}\\
=&\sum\limits_{l=1}^{r-1}\binom{r+s-l-2}{r-l}\binom{\alpha_1}{a}\prod\limits_{j=l+2}^{r+s}\delta_{\alpha_j,0}
+\binom{\alpha_1}{a}\prod\limits_{j=r+2}^{r+s}\delta_{\alpha_j,0}.
\end{align*}
Hence we get
\begin{align*}
c(a,r;b,s)=&\sum\limits_{l=1}^{r-1}\binom{\alpha_1}{a}\binom{r+s-l-1}{r-l}\prod\limits_{j=l+2}^{r+s}\delta_{\alpha_j,0}
+\binom{\alpha_1}{a}\prod\limits_{j=r+2}^{r+s}\delta_{\alpha_j,0}\\
=&\sum\limits_{l=1}^{r}\binom{\alpha_1}{a}\binom{r+s-l-1}{r-l}\prod\limits_{j=l+2}^{r+s}\delta_{\alpha_j,0}.
\end{align*}
Finally, we have the shuffle product formula.

\begin{prop}\label{Prop:Shuffle-Res-1-1}
For any nonnegative integers $a,b$ and any positive integers $r,s$, we have\begin{align}
&x^ay^r\shuffle x^by^s=\sum\limits_{\alpha_1+\cdots+\alpha_{r+s}=a+b\atop \alpha_1,\ldots,\alpha_{r+s}\geqslant 0}\left\{\sum\limits_{l=1}^{r}\binom{\alpha_1}{a}\binom{r+s-l-1}{r-l}\prod\limits_{j=l+2}^{r+s}\delta_{\alpha_j,0}\right.\nonumber\\
&\qquad\qquad\left.+\sum\limits_{l=1}^{s}\binom{\alpha_1}{b}\binom{r+s-l-1}{s-l}\prod\limits_{j=l+2}^{r+s}\delta_{\alpha_j,0}\right\}x^{\alpha_1}y\cdots x^{\alpha_{r+s}}y.
\label{Eq:Shuffle-Res-1-1}
\end{align}
\end{prop}

If $a,b\geqslant 1$, after applying the map $\zeta$ to \eqref{Eq:Shuffle-Res-1-1}, we get the shuffle product formula \cite[Equation (1.1)]{Eie-Wei} of multiple zeta values of the form $\zeta(m,\{1\}^n)$.

It seems that it is tedious and not natural to deduce the shuffle product formula \eqref{Eq:Shuffle-Res-1-1} from Theorem \ref{Thm:Shuffle-General}. In fact, we find that we can get \eqref{Eq:Shuffle-Res-1-1} immediately from the combinatorial description of shuffle product. We want to compute $x^ay^r\shuffle x^by^s$, and we write this product as $x^ay_1^r\shuffle x^by_2^s$. There are two possibilities for the positions of $y$'s:
\begin{description}
  \item[(i)] $\underbrace{y_1\cdots y_1}_l y_2\begin{array}{l}
\underbrace{y_{1}\cdots y_{1}}_{r-l}\\
\underbrace{y_2 \cdots y_{2}}_{s-1}
\end{array}$, \qquad ($1\leqslant l\leqslant r$);
  \item[(ii)] $\underbrace{y_2\cdots y_{2}}_{l}y_1\begin{array}{l}
\underbrace{y_1\cdots y_{1}}_{r-1}\\
\underbrace{y_{2}\cdots y_2}_{s-l}
\end{array}$,\qquad ($1\leqslant l\leqslant s$).
\end{description}
Here and below $\begin{array}{l}
\underbrace{y_{1}\cdots y_{1}}_{m}\\
\underbrace{y_2 \cdots y_{2}}_{n}
\end{array}$ means that $y_{1}^{m}\shuffle y_2^{n}$. For case (i), in $x^{\alpha_1}y$ there are $a$ $x$'s from $x^a$, and $\alpha_1-a$ $x$'s from $x^b$. For $2\leqslant j\leqslant l+1$, all $x$'s  in $x^{\alpha_j}y$ are from $x^b$. And for $l+2\leqslant j\leqslant r+s$, it must hold $\alpha_j=0$. Hence we know that the contribution of case (i) to the coefficient of $x^{\alpha_1}y\cdots x^{\alpha_{r+s}}y$ is
$$\sum\limits_{l=1}^r\binom{\alpha_1}{a}\binom{r+s-l-1}{r-l}\prod\limits_{j=l+2}^{r+s}\delta_{\alpha_j,0}.$$
By the symmetry of $y_1$'s and $y_2$'s, we get that the contribution of case (ii) to the coefficient of $x^{\alpha_1}y\cdots x^{\alpha_{r+s}}y$ is
$$\sum\limits_{l=1}^s\binom{\alpha_1}{b}\binom{r+s-l-1}{s-l}\prod\limits_{j=l+2}^{r+s}\delta_{\alpha_j,0}.$$
Then we reproved the shuffle product formula \eqref{Eq:Shuffle-Res-1-1} in a simple way.

In \cite[Theorem 1.1]{Lei-Guo-Ma}, P. Lei, L. Guo and B. Ma also studied the shuffle product of two multiple zeta values of the form $\zeta(m,\{1\}^n)$. The formula they obtained looks different from that deduced from \eqref{Eq:Shuffle-Res-1-1}. We show in Appendix \ref{AppSec:Proof-1-1} in fact these two formulas essentially coincide with each other.

\subsection{The formula of $x^ay^r\shuffle x^{b_1}y^{s_1}x^{b_2}y^{s_2}$}

In this subsection, we compute the shuffle product $x^ay^r\shuffle x^{b_1}y^{s_1}x^{b_2}y^{s_2}$, where $a,b_1,b_2$ are nonnegative integers and $r,s_1,s_2$ are positive integers. To distinguish these $y$'s, we write
$$x^ay^r=x^ay_1^r,\quad x^{b_1}y^{s_1}x^{b_2}y^{s_2}=x^{b_1}y_2^{s_1}x^{b_2}y_2^{s_2}.$$
There are four cases of the order of $y$'s:
\begin{description}
  \item[(i)] $\underbrace{y_1 \cdots y_{1}}_{r_1}y_2\begin{array}{l}
\underbrace{y_{1}\cdots y_{1}}_{r_2}\\
\underbrace{y_2\cdots y_{2}}_{s_1-2}
\end{array}y_{2}\underbrace{y_{1}\cdots y_{1}}_{r_3}y_2\begin{array}{l}
\underbrace{y_{1} \cdots y_{1}}_{r_4}\\
\underbrace{y_2\cdots y_2}_{s_2-1}
\end{array}$,\\
where $r_1+r_2+r_3+r_4=r$ with $r_1\geqslant 1$ and $r_2,r_3,r_4\geqslant 0$;
  \item[(ii)] $\underbrace{y_2\cdots y_2}_l y_1\begin{array}{l}
\underbrace{y_1 \cdots y_{1}}_{r_1-1}\\
\underbrace{y_{2}\cdots y_{2}}_{s_1-l-1}
\end{array}y_{2}\underbrace{y_{1}\cdots y_{1}}_{r_2}y_2\begin{array}{l}
\underbrace{y_{1}\cdots y_1}_{r_3}\\
\underbrace{y_2\cdots y_{2}}_{s_2-1}
\end{array}$,\\
where $r_1+r_2+r_3=r$ with $r_1\geqslant 1$, $r_2,r_3\geqslant 0$ and $1\leqslant l\leqslant s_1-1$;
  \item[(iii)] $\underbrace{y_2\cdots y_{2}}_{s_1}\underbrace{y_1\cdots y_{1}}_{r_1}y_2\begin{array}{l}
\underbrace{y_{1}\cdots y_1}_{r_2}\\
\underbrace{y_2\cdots y_2}_{s_2-1}
\end{array}$,\\
where $r_1+r_2=r$ with $r_1\geqslant 1$ and $r_2\geqslant 0$;
  \item[(iv)] $\underbrace{y_2\cdots y_{2}}_{s_1}\underbrace{y_2\cdots y_2}_l y_1\begin{array}{l}
\underbrace{y_{1}\cdots y_1}_{r-1}\\
\underbrace{y_{2}\cdots y_2}_{s_2-l}
\end{array}$,\\
where $1\leqslant l\leqslant s_2$.
\end{description}

For each case, we consider the place and the number of $x$'s to get $x^{\alpha_i}y$. Then we can write the contribution of each case to the coefficient of $x^{\alpha_1}y\cdots x^{\alpha_{r+s_1+s_2}}y$ in the product $x^ay^r\shuffle x^{b_1}y^{s_1}x^{b_2}y^{s_2}$, and get the following shuffle product formula.

\begin{prop}
For any nonnegative integers $a,b_1,b_2$ and any positive integers $r,s_1,s_2$, we have
\begin{align}
&x^{a}y^r\shuffle x^{b_1}y^{s_1}x^{b_2}y^{s_2}=\sum\limits_{\alpha_1+\cdots+\alpha_{r+s_1+s_2}=a+b_1+b_2\atop \alpha_1,\ldots,\alpha_{r+s_1+s_2}\geqslant 0}\left\{\sum\limits_{r_1+r_2+r_3+r_4=r\atop r_1\geqslant 1,r_2,r_3,r_4\geqslant 0}\binom{\alpha_1}{a}\binom{r_2+s_1-2}{r_2}\right.\nonumber\\
&\quad \times\binom{r_4+s_2-1}{r_4}\delta_{\alpha_1+\cdots+\alpha_{r_1+1},a+b_1}\prod\limits_{i=r_1+2}^{r_1+r_2+s_1}\delta_{\alpha_i,0}
\prod\limits_{i=r_1+r_2+r_3+s_1+2}^{r+s_1+s_2}\delta_{\alpha_i,0}\nonumber\\
&+\sum\limits_{{r_1+r_2+r_3=r\atop r_1\geqslant 1,r_2,r_3\geqslant 0}\atop 1\leqslant l\leqslant s_1-1}\binom{\alpha_1}{b_1}\binom{r_1+s_1-l-2}{r_1-1}\binom{r_3+s_2-1}{r_3}\delta_{\alpha_1+\cdots+\alpha_{l+1},a+b_1}\prod\limits_{i=l+2}^{r_1+s_1}\delta_{\alpha_i,0}
\nonumber\\
&\quad\times\prod\limits_{i=r_1+r_2+s_1+2}^{r+s_1+s_2}\delta_{\alpha_i,0}+\sum\limits_{r_1+r_2=r\atop r_1\geqslant 1,r_2\geqslant 0}\binom{\alpha_1}{b_1}\binom{\alpha_{s_1+1}}{a+b_1-\sum\limits_{i=1}^{s_1}\alpha_i}\binom{r_2+s_2-1}{r_2}\nonumber\\
&\quad\left.\times\prod\limits_{i=r_1+s_1+2}^{r+s_1+s_2}\delta_{\alpha_i,0}+\sum\limits_{l=1}^{s_2}\binom{\alpha_1}{b_1}\binom{\alpha_{s_1+1}}{b_2}\binom{r+s_2-l-1}{r-1}
\prod\limits_{i=s_1+l+2}^{r+s_1+s_2}\delta_{\alpha_i,0}\right\}\nonumber\\
&\quad\qquad \times x^{\alpha_1}y\cdots x^{\alpha_{r+s_1+s_2}}y.
\label{Eq:Shuffle-Res-1-2}
\end{align}
\end{prop}

When $a,b_1\geqslant 1$, applying the map $\zeta$ to the formula \eqref{Eq:Shuffle-Res-1-2}, we get the formula for the product of the form $\zeta(m,\{1\}^n)\zeta(m_1,\{1\}^{n_1},m_2,\{1\}^{n_2})$. The formula obtained here looks different from that in \cite[Theorem 1.3]{Lei-Guo-Ma}. We show in Appendix \ref{AppSec:Proof-1-2} that these two formulas are essentially identical.

\subsection{The formula of $x^{a_1}y^{r_1}x^{a_2}y^{r_2}\shuffle x^{b_1}y^{s_1}x^{b_2}y^{s_2}$}

In this subsection, we study the restricted shuffle product formula of the product $x^{a_1}y^{r_1}x^{a_2}y^{r_2}\shuffle x^{b_1}y^{s_1}x^{b_2}y^{s_2}$, where $a_1,a_2,b_1,b_2$ are nonnegative integers and $r_1,r_2,s_1,s_2$ are positive integers. As before, we write the $y$ in $x^{a_1}y^{r_1}x^{a_2}y^{r_2}$ as $y_1$ and $y$ in $x^{b_1}y^{s_1}x^{b_2}y^{s_2}$ as $y_2$. When shuffling $y_1$'s and $y_2$'s, by symmetric, we only need consider the cases beginning with $y_1$. Then there are ten cases:

\begin{description}
  \item[(i)] $\underbrace{y_1\cdots y_1}_{r_1}\underbrace{y_1\cdots y_1}_{l_1}y_2\begin{array}{l}
  \underbrace{y_1\cdots y_1}_{l_2}\\
  \underbrace{y_2\cdots y_2}_{s_1-2}
  \end{array}y_2\underbrace{y_1\cdots y_1}_{l_3}y_2\begin{array}{l}
  \underbrace{y_1\cdots y_1}_{l_4}\\
  \underbrace{y_2\cdots y_2}_{s_2-1}
  \end{array}$,\\
where $l_1+l_2+l_3+l_4=r_2$ with $l_1\geqslant 1$, $l_2,l_3,l_4\geqslant 0$;
  \item[(ii)] $\underbrace{y_1\cdots y_1}_{r_1}\underbrace{y_2\cdots y_2}_ky_1\begin{array}{l}
  \underbrace{y_1\cdots y_1}_{l_1-1}\\
  \underbrace{y_2\cdots y_2}_{s_1-k-1}
  \end{array}y_2\underbrace{y_1\cdots y_1}_{l_2}y_2\begin{array}{l}
  \underbrace{y_1\cdots y_1}_{l_3}\\
  \underbrace{y_2\cdots y_2}_{s_2-1}
  \end{array}$,\\
  where $l_1+l_2+l_3=r_2$ with $l_1\geqslant 1$, $l_2,l_3\geqslant 0$ and $1\leqslant k\leqslant s_1-1$;
  \item[(iii)] $\underbrace{y_1\cdots y_1}_{r_1}\underbrace{y_2\cdots y_2}_{s_1}\underbrace{y_1\cdots y_1}_{l_1}y_2\begin{array}{l}
  \underbrace{y_1\cdots y_1}_{l_2}\\
  \underbrace{y_2\cdots y_2}_{s_2-1}
  \end{array}$,\\
  where $l_1+l_2=r_2$ with $l_1\geqslant 1$ and $l_2\geqslant 0$;
  \item[(iv)] $\underbrace{y_1\cdots y_1}_{r_1}\underbrace{y_2\cdots y_2}_{s_1}\underbrace{y_2\cdots y_2}_ky_1\begin{array}{l}
  \underbrace{y_1\cdots y_1}_{r_2-1}\\
  \underbrace{y_2\cdots y_2}_{s_2-k}
  \end{array}$,\\
  where $1\leqslant k\leqslant s_2$;
  \item[(v)] $\underbrace{y_1\cdots y_1}_ly_2\begin{array}{l}
  \underbrace{y_1\cdots y_1}_{r_1-l-1}\\
  \underbrace{y_2\cdots y_2}_{k_1-1}
  \end{array}y_1\underbrace{y_2\cdots y_2}_{k_2}y_1\begin{array}{l}
  \underbrace{y_1\cdots y_1}_{l_1-1}\\
  \underbrace{y_2\cdots y_2}_{k_3-1}
  \end{array}y_2\underbrace{y_1\cdots y_1}_{l_2}y_2\begin{array}{l}
  \underbrace{y_1\cdots y_1}_{l_3}\\
  \underbrace{y_2\cdots y_2}_{s_2-1}
  \end{array}$,\\
  where $1\leqslant l\leqslant r_1-1$, $k_1+k_2+k_3=s_1$ with $k_1,k_3\geqslant 1$, $k_2\geqslant 0$, and $l_1+l_2+l_3=r_2$ with $l_1\geqslant 1$, $l_2,l_3\geqslant 0$;
  \item[(vi)] $\underbrace{y_1\cdots y_1}_{l_1}y_2\begin{array}{l}
  \underbrace{y_1\cdots y_1}_{r_1-l_1-1}\\
  \underbrace{y_2\cdots y_2}_{k-1}
  \end{array}y_1\underbrace{y_2\cdots y_2}_{s_1-k}\underbrace{y_1\cdots y_1}_{l_2}y_2\begin{array}{l}
  \underbrace{y_1\cdots y_1}_{r_2-l_2}\\
  \underbrace{y_2\cdots y_2}_{s_2-1}
  \end{array}$,\\
  where $1\leqslant l_1\leqslant r_1-1$, $1\leqslant l_2\leqslant r_2$ and $1\leqslant k\leqslant s_1-1$;
  \item[(vii)] $\underbrace{y_1\cdots y_1}_{l}y_2\begin{array}{l}
  \underbrace{y_1\cdots y_1}_{r_1-l-1}\\
  \underbrace{y_2\cdots y_2}_{k_1-1}
  \end{array}y_1\underbrace{y_2\cdots y_2}_{s_1-k_1}\underbrace{y_2\cdots y_2}_{k_2}y_1\begin{array}{l}
  \underbrace{y_1\cdots y_1}_{r_2-1}\\
  \underbrace{y_2\cdots y_2}_{s_2-k_2}
  \end{array}$,\\
  where $1\leqslant l\leqslant r_1-1$, $1\leqslant k_1\leqslant s_1-1$ and $1\leqslant k_2\leqslant s_2$;
  \item[(viii)] $\underbrace{y_1\cdots y_1}_{l_1}y_2\begin{array}{l}
  \underbrace{y_1\cdots y_1}_{l_2}\\
  \underbrace{y_2\cdots y_2}_{s_1-2}
  \end{array}y_2\underbrace{y_1\cdots y_1}_{l_3}\underbrace{y_1\cdots y_1}_{l}y_2\begin{array}{l}
  \underbrace{y_1\cdots y_1}_{r_2-l}\\
  \underbrace{y_2\cdots y_2}_{s_2-1}
  \end{array}$,\\
  where $l_1+l_2+l_3=r_1$ with $l_1,l_3\geqslant 1$, $l_2\geqslant 0$, and $1\leqslant l\leqslant r_2$;
  \item[(ix)] $\underbrace{y_1\cdots y_1}_{l_1}y_2\begin{array}{l}
  \underbrace{y_1\cdots y_1}_{l_2}\\
  \underbrace{y_2\cdots y_2}_{s_1-2}
  \end{array}y_2\underbrace{y_1\cdots y_1}_{l_3}\underbrace{y_2\cdots y_2}_ky_1\begin{array}{l}
  \underbrace{y_1\cdots y_1}_{r_2-1}\\
  \underbrace{y_2\cdots y_2}_{s_2-k}
  \end{array}$,\\
  where $l_1+l_2+l_3=r_1$ with $l_1,l_3\geqslant 1$, $l_2\geqslant 0$, and $1\leqslant k\leqslant s_2$;
  \item[(x)] $\underbrace{y_1\cdots y_1}_{l_1}y_2\begin{array}{l}
  \underbrace{y_1\cdots y_1}_{l_2}\\
  \underbrace{y_2\cdots y_2}_{s_1-2}
  \end{array}y_2\underbrace{y_1\cdots y_1}_{l_3}y_2\begin{array}{l}
  \underbrace{y_1\cdots y_1}_{l_4-1}\\
  \underbrace{y_2\cdots y_2}_{k_1-1}
  \end{array}y_1\underbrace{y_2\cdots y_2}_{k_2}y_1\begin{array}{l}
  \underbrace{y_1\cdots y_1}_{r_2-1}\\
  \underbrace{y_2\cdots y_2}_{k_3}
  \end{array}$,\\
  where $l_1+l_2+l_3+l_4=r_1$ with $l_1,l_4\geqslant 1$, $l_2,l_3\geqslant 0$, and $k_1+k_2+k_3=s_2$ with $k_1\geqslant 1$, $k_2,k_3\geqslant 0$.
\end{description}

Then consider the possible positions of $x$'s, we get the following shuffle product formula.

\begin{prop}
For any nonnegative integers $a_1,a_2,b_1,b_2$ and any positive integers $r_1,r_2$, $s_1,s_2$, we have
\begin{align}
&x^{a_1}y^{r_1}x^{a_2}y^{r_2}\shuffle x^{b_1}y^{s_1}x^{b_2}y^{s_2}=\sum\limits_{{\alpha_1+\cdots+\alpha_{r_1+r_2+s_1+s_2}\atop=a_1+a_2+b_1+b_2}\atop \alpha_1,\ldots,\alpha_{r_1+r_2+s_1+s_2}\geqslant 0}[c(a_1,a_2,r_1,r_2;b_1,b_2,s_1,s_2)\nonumber\\
&\qquad +c(b_1,b_2,s_1,s_2;a_1,a_2,r_1,r_2)]x^{\alpha_1}y\cdots x^{\alpha_{r_1+r_2+s_1+s_2}}y,
\label{Eq:Shuffle-Res-2-2}
\end{align}
where
\begin{align*}
&c(a_1,a_2,r_1,r_2;b_1,b_2,s_1,s_2)=\sum\limits_{l_1+l_2+l_3+l_4=r_2\atop l_1\geqslant 1,l_2,l_3,l_4\geqslant 0}\binom{\alpha_1}{a_1}\binom{\alpha_{r_1+1}}{a_2}\binom{l_2+s_1-2}{l_2}\\
&\;\times \binom{l_4+s_2-1}{l_4}\delta_{\alpha_1+\cdots+\alpha_{r_1+l_1+1},a_1+a_2+b_1}
\prod\limits_{i=r_1+l_1+2}^{r_1+l_1+l_2+s_1}\delta_{\alpha_i,0}\prod\limits_{i=r_1+l_1+l_2+l_3+s_1+2}^{r_1+r_2+s_1+s_2}\delta_{\alpha_i,0}\\
&\;+\sum\limits_{{l_1+l_2+l_3=r_2\atop l_1\geqslant 1,l_2,l_3\geqslant 0}\atop 1\leqslant k\leqslant s_1-1}\binom{\alpha_1}{a_1}\binom{\alpha_{r_1+1}}{a_1+b_1-\sum\limits_{i=1}^{r_1}\alpha_i}\binom{l_1+s_1-k-2}{l_1-1}\binom{l_3+s_2-1}{l_3}\\
&\;\times\delta_{\alpha_1+\cdots+\alpha_{r_1+k+1},a_1+a_2+b_1}\prod\limits_{i=r_1+k+2}^{r_1+l_1+s_1}\delta_{\alpha_i,0}\prod\limits_{i=r_1+l_1+l_2+s_1+2}^{r_1+r_2+s_1+s_2}
\delta_{\alpha_i,0}\\
&\;+\sum\limits_{l_1+l_2=r_2\atop l_1\geqslant 1,l_2\geqslant 0}\binom{\alpha_1}{a_1}\binom{\alpha_{r_1+1}}{a_1+b_1-\sum\limits_{i=1}^{r_1}\alpha_i}\binom{\alpha_{r_1+s_1+1}}{a_1+a_2+b_1-\sum\limits_{i=1}^{r_1+s_1}\alpha_i}
\binom{l_2+s_2-1}{l_2}\\
&\;\times\prod\limits_{i=r_1+l_1+s_1+2}^{r_1+r_2+s_1+s_2}\delta_{\alpha_i,0}+\sum\limits_{k=1}^{s_2}\binom{\alpha_1}{a_1}
\binom{\alpha_{r_1+1}}{a_1+b_1-\sum\limits_{i=1}^{r_1}\alpha_i}\binom{\alpha_{r_1+s_1+1}}{b_2}\\
&\;\times\binom{r_2+s_2-k-1}{r_2-1}\prod\limits_{i=r_1+s_1+k+2}^{r_1+r_2+s_1+s_2}\delta_{\alpha_i,0}+\sum\limits_{{k_1+k_2+k_3=s_1\atop l_1+l_2+l_3=r_2}\atop {k_1,k_3,l_1\geqslant 1,k_2,l_2,l_3\geqslant 0\atop 1\leqslant l\leqslant r_1-1}}\binom{\alpha_1}{a_1}\\
&\;\times \binom{r_1+k_1-l-2}{k_1-1}\binom{\alpha_{r_1+k_1+k_2+1}}{a_2-\sum\limits_{i=r_1+k_1+1}^{r_1+k_1+k_2}\alpha_i}\binom{k_3+l_1-2}{l_1-1}
\binom{l_3+s_2-1}{l_3}\\
&\;\times \delta_{\alpha_1+\cdots+\alpha_{l+1},a_1+b_1}\prod\limits_{i=l+2}^{r_1+k_1}\delta_{\alpha_i,0}\prod\limits_{i=r_1+k_1+k_2+2}^{r_1+k_1+k_2+k_3+l_1}\delta_{\alpha_i,0}
\prod\limits_{i=r_1+l_1+l_2+s_1+2}^{r_1+r_2+s_1+s_2}\delta_{\alpha_i,0}\\
&\;+\sum\limits_{{1\leqslant l_1\leqslant r_1-1\atop 1\leqslant l_2\leqslant r_2}\atop 1\leqslant k\leqslant s_1-1}\binom{\alpha_1}{a_1}\binom{k+r_1-l_1-2}{k-1}\binom{\alpha_{r_1+s_1+1}}{a_2-\sum\limits_{i=k+r_1+1}^{r_1+s_1}\alpha_i}\binom{r_2+s_2-l_2-1}{s_2-1}\\
&\times\delta_{\alpha_1+\cdots+\alpha_{l_1+1},a_1+b_1}\prod\limits_{i=l_1+2}^{k+r_1}\delta_{\alpha_i,0}\prod\limits_{i=r_1+l_2+s_1+2}^{r_1+r_2+s_1+s_2}\delta_{\alpha_i,0}
+\sum\limits_{{1\leqslant l\leqslant r_1-1\atop 1\leqslant k_1\leqslant s_1-1}\atop 1\leqslant k_2\leqslant s_2}\binom{\alpha_1}{a_1}\\
&\;\times\binom{k_1+r_1-l-2}{k_1-1}\binom{\alpha_{r_1+s_1+1}}{b_2}\binom{r_2+s_2-k_2-1}{r_2-1}\delta_{\alpha_1+\cdots+\alpha_{l+1},a_1+b_1}\\
&\;\times\prod\limits_{i=l+2}^{k_1+r_1}\delta_{\alpha_i,0}\prod\limits_{i=k_2+r_1+s_1+2}^{r_1+r_2+s_1+s_2}\delta_{\alpha_i,0}+\sum\limits_{{l_1+l_2+l_3=r_1\atop l_1,l_3\geqslant 1,l_2\geqslant 0}\atop 1\leqslant l\leqslant r_2}\binom{\alpha_1}{a_1}\binom{l_2+s_1-2}{l_2}\binom{\alpha_{r_1+s_1+1}}{a_2}\\
&\;\times \binom{r_2+s_2-l-1}{s_2-1}\delta_{\alpha_1+\cdots+\alpha_{l_1+1},a_1+b_1}
\prod\limits_{i=l_1+2}^{l_1+l_2+s_1}\delta_{\alpha_i,0}\prod\limits_{i=r_1+s_1+l+2}^{r_1+r_2+s_1+s_2}\delta_{\alpha_i,0}\\
&\;+\sum\limits_{{l_1+l_2+l_3=r_1\atop l_1,l_3\geqslant 1,l_2\geqslant 0}\atop 1\leqslant k\leqslant s_2}\binom{\alpha_1}{a_1}\binom{l_2+s_1-2}{l_2}\binom{\alpha_{r_1+s_1+1}}{b_2-\sum\limits_{i=l_1+l_2+s_1+1}^{r_1+s_1}\alpha_i}\binom{r_2+s_2-k-1}{r_2-1}\\
&\;\times \delta_{\alpha_1+\cdots+\alpha_{l_1+1},a_1+b_1}
\prod\limits_{i=l_1+2}^{l_1+l_2+s_1}\delta_{\alpha_i,0}\prod\limits_{i=k+r_1+s_1+2}^{r_1+r_2+s_1+s_2}\delta_{\alpha_i,0}+\sum\limits_{{l_1+l_2+l_3+l_4=r_1\atop k_1+k_2+k_3=s_2}\atop l_1,l_4,k_1\geqslant 1,l_2,l_3,k_2,k_3\geqslant 0}\binom{\alpha_1}{a_1}\\
&\;\times \binom{l_2+s_1-2}{l_2}
\binom{\alpha_{l_1+l_2+l_3+s_1+1}}{b_2-\sum\limits_{i=l_1+l_2+s_1+1}^{l_1+l_2+l_3+s_1}\alpha_i}\binom{k_1+l_4-2}{l_4-1}\binom{k_3+r_2-1}{k_3}\\
&\;\times\delta_{\alpha_1+\cdots+\alpha_{l_1+1},a_1+b_1}\prod\limits_{i=l_1+2}^{l_1+l_2+s_1}\delta_{\alpha_i,0}\prod\limits_{i=l_1+l_2+l_3+s_1+2}^{k_1+r_1+s_1}\delta_{\alpha_i,0}
\prod\limits_{i=k_1+k_2+r_1+s_1+2}^{r_1+r_2+s_1+s_2}\delta_{\alpha_i,0}.
\end{align*}
\end{prop}


\section{Shuffle products of several multiple zeta values}\label{Sec:Shu-MZV-n}

There is another approach to get \eqref{Eq:Shuffle-Res-1-1} of the product $x^ay^r\shuffle x^by^s$. We want to find the coefficient of $x^{\alpha_1}y\cdots x^{\alpha_{r+s}}y$ in the product $x^ay^r\shuffle x^by^s$ for $\alpha_1+\cdots+\alpha_{r+s}=a+b$ with $\alpha_i\geqslant 0$. We write this product as $x_1^ay_1^r\shuffle x_2^by_2^s$, and consider $y$'s first and $x$'s second. To get $r+s$ $y$'s, by symmetry, we only need to consider the case
$$\fbox{$y_1$} \underbrace{y\cdots y}_{l-1} \fbox{$y_2$}\underbrace{y\cdots y}_{r+s-l-1},$$
where $1\leqslant l\leqslant r$ and $\fbox{$y_1$}$ denotes the first $y_1$. The $y$'s between $\fbox{$y_1$}$ and $\fbox{$y_2$}$ can only be $y_1$, while the $y$'s after $\fbox{$y_2$}$ are either $y_1$ or $y_2$. The position of $\fbox{$y_2$}$ is fixed, and there are $r+s-l-1$ positions to place other $y_2$'s. After all $y_2$'s are placed, there is only one way to place $y_1$'s. Hence in this case, there are $\binom{r+s-l-1}{s-1}$ possibilities to get $y^{r+s}$. Now we come to place $x$'s. All $x_1$'s are placed before $\fbox{$y_1$}$, and all $x_2$'s can be placed free before $\fbox{$y_2$}$. Hence there are no $x$'s before the $y$'s which are after $\fbox{$y_2$}$. In other words we have $\alpha_i=0$ for $l+2\leqslant i\leqslant r+s$. There are $\alpha_1$ positions to place $x_1$'s, and after all $x_1$'s are placed, there is only one way to place $x_2$'s. Hence we get $\binom{\alpha_1}{a}$ possibilities. Finally, we get the contribution of this case to the coefficient of $x^{\alpha_1}y\cdots x^{\alpha_{r+s}}y$ is
$$\sum\limits_{l=1}^r\binom{\alpha_1}{a}\binom{r+s-l-1}{s-1}\prod\limits_{i=l+2}^{r+s}\delta_{\alpha_i,0}.$$
Then we get \eqref{Eq:Shuffle-Res-1-1}.

Generalizing the arguments above, we can get the formula of the shuffle products $x^{a_1}y^{r_1}\shuffle x^{a_2}y^{r_2}\shuffle\cdots\shuffle x^{a_n}y^{r_n}$.

\begin{prop}
Let $n,r_1,\ldots,r_n$ be positive integers and let $a_1,\ldots,a_n$ be nonnegative integers. Then we have
\begin{align}
&x^{a_1}y^{r_1}\shuffle x^{a_2}y^{r_2}\shuffle\cdots\shuffle x^{a_n}y^{r_n}\nonumber\\
=&\sum\limits_{\alpha_1+\cdots+\alpha_{r_1+\cdots+r_n}=a_1+\cdots+a_n\atop \alpha_i\geqslant 0}c_{\alpha_1,\ldots,\alpha_{r_1+\cdots+r_n}}x^{\alpha_1}y\cdots x^{\alpha_{r_1+\cdots+r_n}}y,
\label{Eq:Shuffle-Res-nProd}
\end{align}
with the coefficient
\begin{align*}
&c_{\alpha_1,\ldots,\alpha_{r_1+\cdots+r_n}}=\sum\limits_{l_1+\cdots+l_n=r_1+\cdots+r_n\atop l_i\geqslant 1}\sum\limits_{\sigma\in\mathfrak{S}_n}\sigma_r\left\{\prod\limits_{j=2}^n\begin{pmatrix}
{\sum\limits_{i=j}^nl_i-\sum\limits_{i=j+1}^nr_i-1}\\
{r_j-1}
\end{pmatrix}\right\}\\
&\quad\;\times\sigma_a\left\{\prod\limits_{j=1}^{n-1}\begin{pmatrix}
{\sum\limits_{i=1}^{L_{j-1}+1}\alpha_i-\sum\limits_{i=1}^{j-1}a_i}\\
{a_j}
\end{pmatrix}\right\}\prod\limits_{j=L_{n-1}+2}^{r_1+\cdots+r_n}
\delta_{\alpha_i,0},
\end{align*}
where $L_0=0$ and $L_j=l_1+\cdots+l_j$ for $1\leqslant j\leqslant n$, and $\sigma_r$, $\sigma_a$ are induced permutations of $\sigma\in\mathfrak{S}_n$ on the set $\{r_1,\ldots,r_n\}$, $\{a_1,\ldots,a_n\}$, respectively.
\end{prop}

\proof We write this product as $x_1^{a_1}y_1^{r_1}\shuffle x_2^{a_2}y_2^{r_2}\shuffle\cdots\shuffle x_n^{a_n}y_n^{r_n}$, and consider $y$'s first and $x$'s second to get $x^{\alpha_1}y\cdots x^{\alpha_{r_1+\cdots+r_n}}y$ from the combinatorial description of shuffle product. By symmetry, we only need to consider the case
$$\fbox{$y_1$}\underbrace{y\cdots y}_{l_1-1}\fbox{$y_2$}\underbrace{y\cdots y}_{l_2-1}\cdots\fbox{$y_{n-1}$}\underbrace{y\cdots y}_{l_{n-1}-1}\fbox{$y_n$}\underbrace{y\cdots y}_{l_n-1},$$
where $l_1+\cdots+l_n=r_1+\cdots+r_n$ with $l_i\geqslant 1$.

We compute the contribution of this case to the coefficient $c_{\alpha_1,\ldots,\alpha_{r_1+\cdots+r_n}}$. The position of $\fbox{$y_n$}$ is fixed, and there are $l_n-1$ positions to place other $y_n$'s. Hence the possibility is $\binom{l_n-1}{r_n-1}$. For $y_{n-1}$, the position of $\fbox{$y_{n-1}$}$ is fixed, and there are $l_{n-1}+l_n-1$ positions to place other $y_{n-1}$'s except that $r_n$ positions are occupied by $y_n$'s. Then the possibility is $\binom{l_{n-1}+l_n-r_n-1}{r_{n-1}-1}$. Repeating this arguments, for $y_2$'s, the position of $\fbox{$y_2$}$ is fixed, and there are $l_2+\cdots+l_n-1$ positions to place other $y_2$'s except that $r_3+\cdots+r_n$ positions are occupied by $y_3$'s, $\ldots$, $y_n$'s. Hence the possibility is $\binom{l_2+\cdots+l_n-r_3-\cdots-r_n-1}{r_2-1}$. After all $y_2$'s, $\ldots$, $y_n$'s are placed, there is only one way to place $y_1$'s. Then we get the product
$$\prod\limits_{j=2}^n\begin{pmatrix}
{\sum\limits_{i=j}^nl_i-\sum\limits_{i=j+1}^nr_i-1}\\
{r_j-1}
\end{pmatrix}.$$

Now we come to place $x$'s. There is $\alpha_1$ positions to place $x_1$. Hence the possibility is $\binom{\alpha_1}{a_1}$. For $x_2$, there are $\alpha_1+\cdots+\alpha_{L_1+1}$ positions to place $x_2$'s except that $a_1$ positions are occupied by $x_1$'s. Hence the possibility is $\binom{\alpha_1+\cdots+\alpha_{L_1+1}-a_1}{a_2}$. Repeating this arguments, for $x_{n-1}$, there are $\alpha_1+\cdots+\alpha_{L_{n-2}+1}$ positions to place $x_{n-1}$, except that $a_1+\cdots+a_{n-2}$ positions are occupied by $x_1$'s, $\ldots$, $x_{n-2}$'s. Then the possibility is $\binom{\alpha_1+\cdots+\alpha_{L_{n-2}+1}-a_1-\cdots-a_{n-2}}{a_{n-1}}$. After all $x_1$'s, $\ldots$, $x_{n-1}$'s are placed, there is only one way to place $x_n$. Hence we get the product
$$\prod\limits_{j=1}^{n-1}\begin{pmatrix}
{\sum\limits_{i=1}^{L_{j-1}+1}\alpha_i-\sum\limits_{i=1}^{j-1}a_i}\\
{a_j}
\end{pmatrix}.$$

Finally, it is easy to see that for $L_{n-1}+2\leqslant j\leqslant r_1+\cdots+r_n$, it must hold $\alpha_i=0$.
\qed

If $r_1=\cdots=r_n=1$, then the formula \eqref{Eq:Shuffle-Res-nProd} becomes
\begin{align}
&x^{a_1}y\shuffle x^{a_2}y\shuffle\cdots\shuffle x^{a_n}y\nonumber\\
=&\sum\limits_{\alpha_1+\cdots+\alpha_{n}=a_1+\cdots+a_n\atop \alpha_i\geqslant 0}\sum\limits_{\sigma\in\mathfrak{S}_n}\sigma_a\left\{\prod\limits_{j=1}^{n-1}\begin{pmatrix}
{\sum\limits_{i=1}^{j}\alpha_i-\sum\limits_{i=1}^{j-1}a_i}\\
{a_j}
\end{pmatrix}\right\}x^{\alpha_1}y\cdots x^{\alpha_{n}}y.
\label{Eq:Shuffle-nProd}
\end{align}
Under the case $a_1,\ldots,a_n\geqslant 1$, after applying the map $\zeta$, we get \cite[Main Theorem]{Eie-Liaw-Ong} from \eqref{Eq:Shuffle-Res-nProd}, and get \cite[Theorem 1.3]{Chung-Eie-Liaw-Ong} from \eqref{Eq:Shuffle-nProd}.


\appendix


\section{The equivalence of two formulas of the product $x^ay^r\shuffle x^by^s$}\label{AppSec:Proof-1-1}

The shuffle product formula of two multiple zeta values of the form $\zeta(m,\{1\}^n)$ in \cite[Theorem 1.1 or Equation (23)]{Lei-Guo-Ma} can be expressed as
\begin{align}
&x^ay^r\shuffle x^by^s=\sum\limits_{{0\leqslant k\leqslant b \atop r_1+r_2=r,r_i\geqslant 0}\atop \alpha_1+\cdots+\alpha_{r_1+1}=b-k,\alpha_i\geqslant 0}\binom{a-1+k}{a-1}\binom{r_2+s-1}{s-1}x^{\alpha_1+a+k}yx^{\alpha_2}y\nonumber\\
&\quad \cdots x^{\alpha_{r_1+1}}y^{r_2+s}+\sum\limits_{{1\leqslant l\leqslant s \atop a_1+a_2=a-1,a_i\geqslant 0}\atop \alpha_1+\cdots+\alpha_{l+1}=a_2,\alpha_i\geqslant 0}\binom{a_1+b-1}{b-1}\binom{r+s-l}{r}\nonumber\\
&\qquad \times x^{\alpha_1+a_1+b}yx^{\alpha_2}y\cdots x^{\alpha_{l}}yx^{\alpha_{l+1}+1}y^{r+s-l},
\label{Eq:Shuffle-Res-1-1-Old}
\end{align}
where $a,b,r,s$ are positive integers. We show that the shuffle product formula \eqref{Eq:Shuffle-Res-1-1-Old} can be deduced from \eqref{Eq:Shuffle-Res-1-1}.

Denote the first sum of the right-hand side of \eqref{Eq:Shuffle-Res-1-1-Old} by $\Sigma_1$ and the second sum by $\Sigma_2$, then we have
\begin{align*}
\Sigma_1=&\sum\limits_{{0\leqslant k\leqslant b \atop 0\leqslant l\leqslant r}\atop {\alpha_1+\cdots+\alpha_{l+1}=a+b,\atop \alpha_1\geqslant a+k,\alpha_i\geqslant 0}}\binom{a-1+k}{a-1}\binom{r+s-l-1}{s-1}x^{\alpha_1}yx^{\alpha_2}y\cdots x^{\alpha_{l+1}}y^{r+s-l}\\
=&\sum\limits_{{\alpha_1+\cdots+\alpha_{l+1}=a+b}\atop {0\leqslant l\leqslant r,\alpha_i\geqslant 0}}\left(\sum\limits_{k=0}^{\alpha_1-a}\binom{a-1+k}{a-1}\right)\binom{r+s-l-1}{s-1}x^{\alpha_1}yx^{\alpha_2}y\cdots x^{\alpha_{l+1}}y^{r+s-l}\\
=&\sum\limits_{{\alpha_1+\cdots+\alpha_{l+1}=a+b}\atop {0\leqslant l\leqslant r,\alpha_i\geqslant 0}}\binom{\alpha_1}{a}\binom{r+s-l-1}{s-1}x^{\alpha_1}yx^{\alpha_2}y\cdots x^{\alpha_{l+1}}y^{r+s-l},
\end{align*}
and
\begin{align*}
\Sigma_2=&\sum\limits_{{1\leqslant l\leqslant s \atop 0\leqslant k\leqslant a-1}\atop {\alpha_1+\cdots+\alpha_{l+1}=a+b,\atop \alpha_1\geqslant k+b,\alpha_{l+1}\geqslant 1,\alpha_i\geqslant 0}}\binom{k+b-1}{b-1}\binom{r+s-l}{r}x^{\alpha_1}yx^{\alpha_2}y\cdots x^{\alpha_{l}}yx^{\alpha_{l+1}}y^{r+s-l}\\
=&\sum\limits_{{\alpha_1+\cdots+\alpha_{l+1}=a+b,\atop \alpha_{l+1}\geqslant 1,\alpha_i\geqslant 0}\atop 1\leqslant l\leqslant s}\left(\sum\limits_{k=0}^{\alpha_1-b}\binom{k+b-1}{b-1}\right)\binom{r+s-l}{r}x^{\alpha_1}yx^{\alpha_2}y\cdots x^{\alpha_{l+1}}y^{r+s-l}\\
=&\sum\limits_{{\alpha_1+\cdots+\alpha_{l+1}=a+b,\atop \alpha_{l+1}\geqslant 1,\alpha_i\geqslant 0}\atop 1\leqslant l\leqslant s}\binom{\alpha_1}{b}\binom{r+s-l}{r}x^{\alpha_1}yx^{\alpha_2}y\cdots x^{\alpha_{l+1}}y^{r+s-l}.
\end{align*}
Moreover, we have
\begin{align*}
\Sigma_2=&\sum\limits_{\alpha_1+\cdots+\alpha_{l+1}=a+b,\atop \alpha_i\geqslant 0,1\leqslant l\leqslant s}\binom{\alpha_1}{b}\binom{r+s-l}{r}x^{\alpha_1}yx^{\alpha_2}y\cdots x^{\alpha_{l+1}}y^{r+s-l}\\
&-\sum\limits_{\alpha_1+\cdots+\alpha_{l}=a+b,\atop \alpha_i\geqslant 0,1\leqslant l\leqslant s}\binom{\alpha_1}{b}\binom{r+s-l}{r}x^{\alpha_1}yx^{\alpha_2}y\cdots x^{\alpha_{l}}y^{r+s-l+1}\\
=&\sum\limits_{\alpha_1+\cdots+\alpha_{l+1}=a+b,\atop \alpha_i\geqslant 0,1\leqslant l\leqslant s}\binom{\alpha_1}{b}\binom{r+s-l}{r}x^{\alpha_1}yx^{\alpha_2}y\cdots x^{\alpha_{l+1}}y^{r+s-l}\\
&-\sum\limits_{\alpha_1+\cdots+\alpha_{l+1}=a+b,\atop \alpha_i\geqslant 0,0\leqslant l\leqslant s-1}\binom{\alpha_1}{b}\binom{r+s-l-1}{r}x^{\alpha_1}yx^{\alpha_2}y\cdots x^{\alpha_{l+1}}y^{r+s-l}\\
=&\sum\limits_{\alpha_1+\cdots+\alpha_{l+1}=a+b,\atop \alpha_i\geqslant 0,1\leqslant l\leqslant s}\binom{\alpha_1}{b}\binom{r+s-l-1}{r-1}x^{\alpha_1}yx^{\alpha_2}y\cdots x^{\alpha_{l+1}}y^{r+s-l}\\
&\quad -\binom{a+b}{b}\binom{r+s-1}{r}x^{a+b}y^{r+s},
\end{align*}
which implies that
\begin{align*}
\Sigma_1+\Sigma_2=&\sum\limits_{{\alpha_1+\cdots+\alpha_{l+1}=a+b}\atop {1\leqslant l\leqslant r,\alpha_i\geqslant 0}}\binom{\alpha_1}{a}\binom{r+s-l-1}{s-1}x^{\alpha_1}yx^{\alpha_2}y\cdots x^{\alpha_{l+1}}y^{r+s-l}\\
&+\sum\limits_{\alpha_1+\cdots+\alpha_{l+1}=a+b,\atop \alpha_i\geqslant 0,1\leqslant l\leqslant s}\binom{\alpha_1}{b}\binom{r+s-l-1}{r-1}x^{\alpha_1}yx^{\alpha_2}y\cdots x^{\alpha_{l+1}}y^{r+s-l}.
\end{align*}
Then we have shown that the  formula \eqref{Eq:Shuffle-Res-1-1-Old} is a consequence of \eqref{Eq:Shuffle-Res-1-1}.


\section{The equivalence of two formulas of the product $x^ay^r\shuffle x^{b_1}y^{s_1}x^{b_2}y^{s_2}$}\label{AppSec:Proof-1-2}

In \cite[Theorem 1.3 or the last equation]{Lei-Guo-Ma}, the authors showed that for any positive integers $a,b_1,b_2,r,s_1,s_2$, the following shuffle product formula holds
\begin{align}
x^ay^r\shuffle x^{b_1}y^{s_1}x^{b_2}y^{s_2}=\Sigma_1+\Sigma_2+\Sigma_3+\Sigma_4,
\label{Eq:Shuffle-Res-1-2-Old}
\end{align}
where
\begin{align*}
\Sigma_1=&\sum\limits_{{0\leqslant k\leqslant b_1\atop r_1+r_2+r_3+r_4=r,r_i\geqslant 0}\atop {\alpha_1+\cdots+\alpha_{r_1+1}=b_1-k,\alpha_i\geqslant 0\atop \widetilde{\alpha}_1+\cdots+\widetilde{\alpha}_{r_3+1}=b_2-1,\widetilde{\alpha}_j\geqslant 0}}\binom{a-1+k}{a-1}\binom{r_2+s_1-1}{s_1-1}\binom{r_4+s_2-1}{s_2-1}\\
&\; \times x^{\alpha_1+a+k}yx^{\alpha_2}y\cdots x^{\alpha_{r_1}}yx^{\alpha_{r_1+1}}y^{r_2+s_1}x^{\widetilde{\alpha}_1+1}yx^{\widetilde{\alpha}_2}y\cdots x^{\widetilde{\alpha}_{r_3}}yx^{\widetilde{\alpha}_{r_3+1}}y^{r_4+s_2},\\
\Sigma_2=&\sum\limits_{{{1\leqslant l\leqslant s_1\atop a_1+a_2=a-1,a_i\geqslant 0}\atop r_1+r_2+r_3=r,r_i\geqslant 0}\atop{\alpha_1+\cdots+\alpha_{l+1}=a_2,\alpha_i\geqslant 0\atop \widetilde{\alpha}_1+\cdots+\widetilde{\alpha}_{r_2+1}=b_2-1,\widetilde{\alpha}_i\geqslant 0}}\binom{a_1+b_1-1}{b_1-1}\binom{r_1+s_1-l}{s_1-l}\binom{r_3+s_2-1}{s_2-1}\\
&\;\times x^{\alpha_1+a_1+b_1}yx^{\alpha_2}y\cdots x^{\alpha_l}yx^{\alpha_{l+1}+1}y^{r_1+s_1-l}x^{\widetilde{\alpha}_1+1}yx^{\widetilde{\alpha}_2}y\cdots x^{\widetilde{\alpha}_{r_2}}yx^{\widetilde{\alpha}_{r_2+1}}y^{r_3+s_2},\\
\Sigma_3=&\sum\limits_{{{1\leqslant k\leqslant b_2\atop a_1+a_2+a_3=a-1,a_i\geqslant 0}\atop r_1+r_2=r,r_i\geqslant 0}\atop{\alpha_1+\cdots+\alpha_{s_1+1}=a_2,\alpha_i\geqslant 0\atop\widetilde{\alpha}_1+\cdots+\widetilde{\alpha}_{r_1+1}=b_2-k,\widetilde{\alpha}_i\geqslant 0}}\binom{a_1+b_1-1}{b_1-1}\binom{a_3+k-1}{k-1}\binom{r_2+s_2-1}{s_2-1}\\
&\;\times x^{\alpha_1+a_1+b_1}yx^{\alpha_2}y\cdots x^{\alpha_{s_1}}yx^{\alpha_{s_1+1}+\widetilde{\alpha}_1+a_3+k+1}yx^{\widetilde{\alpha}_2}y\cdots x^{\widetilde{\alpha}_{r_1}}yx^{\widetilde{\alpha}_{r_1+1}}y^{r_2+s_2},\\
\Sigma_4=&\sum\limits_{{1\leqslant l\leqslant s_2\atop a_1+a_2+a_3+a_4=a-1,a_i\geqslant 0}\atop{\alpha_1+\cdots+\alpha_{s_1}=a_2,\alpha_i\geqslant 0\atop \widetilde{\alpha}_1+\cdots+\widetilde{\alpha}_{l+1}=a_4,\widetilde{\alpha}_i\geqslant 0}}\binom{a_1+b_1-1}{b_1-1}\binom{a_3+b_2-1}{b_2-1}\binom{r+s_2-l}{r}\\
&\;\times x^{\alpha_1+a_1+b_1}yx^{\alpha_2}y\cdots x^{\alpha_{s_1}}yx^{\widetilde{\alpha}_1+a_3+b_2}yx^{\widetilde{\alpha}_2}y\cdots x^{\widetilde{\alpha}_l}yx^{\widetilde{\alpha}_{l+1}+1}y^{r+s_2-l}.
\end{align*}
We show that the formulas \eqref{Eq:Shuffle-Res-1-2} and \eqref{Eq:Shuffle-Res-1-2-Old} are essentially equivalent.

We have
\begin{align*}
\Sigma_1=&\sum\limits_{{0\leqslant k\leqslant b_1\atop r_1+r_2+r_3+r_4=r}\atop {\alpha_1+\cdots+\alpha_{r_1+1}=a+b_1,\alpha_1\geqslant a+k\atop \widetilde{\alpha}_1+\cdots+\widetilde{\alpha}_{r_3+1}=b_2,\widetilde{\alpha}_1\geqslant 1}}\binom{a-1+k}{a-1}\binom{r_2+s_1-1}{r_2}\binom{r_4+s_2-1}{r_4}\\
&\; \times x^{\alpha_1}y\cdots x^{\alpha_{r_1}}yx^{\alpha_{r_1+1}}y^{r_2+s_1}x^{\widetilde{\alpha}_1}y\cdots x^{\widetilde{\alpha}_{r_3}}yx^{\widetilde{\alpha}_{r_3+1}}y^{r_4+s_2}\\
=&\sum\limits_{{r_1+r_2+r_3+r_4=r}\atop {\alpha_1+\cdots+\alpha_{r_1+1}=a+b_1\atop \widetilde{\alpha}_1+\cdots+\widetilde{\alpha}_{r_3+1}=b_2,\widetilde{\alpha}_1\geqslant 1}}\binom{\alpha_1}{a}\binom{r_2+s_1-1}{r_2}\binom{r_4+s_2-1}{r_4}\\
&\; \times x^{\alpha_1}y\cdots x^{\alpha_{r_1}}yx^{\alpha_{r_1+1}}y^{r_2+s_1}x^{\widetilde{\alpha}_1}y\cdots x^{\widetilde{\alpha}_{r_3}}yx^{\widetilde{\alpha}_{r_3+1}}y^{r_4+s_2}.
\end{align*}
Here and below, all indexes appearing in the summations are nonnegative integers without special statement. Neglecting the condition $\widetilde{\alpha}_1\geqslant 1$, we have
\begin{align*}
\Sigma_1=&\sum\limits_{{r_1+r_2+r_3+r_4=r}\atop {\alpha_1+\cdots+\alpha_{r_1+1}=a+b_1\atop \widetilde{\alpha}_1+\cdots+\widetilde{\alpha}_{r_3+1}=b_2}}\binom{\alpha_1}{a}\binom{r_2+s_1-1}{r_2}\binom{r_4+s_2-1}{r_4}\\
&\; \times x^{\alpha_1}y\cdots x^{\alpha_{r_1}}yx^{\alpha_{r_1+1}}y^{r_2+s_1}x^{\widetilde{\alpha}_1}y\cdots x^{\widetilde{\alpha}_{r_3}}yx^{\widetilde{\alpha}_{r_3+1}}y^{r_4+s_2}\\
&\;-\sum\limits_{{r_1+r_2+r_3+r_4=r}\atop {\alpha_1+\cdots+\alpha_{r_1+1}=a+b_1\atop \widetilde{\alpha}_2+\cdots+\widetilde{\alpha}_{r_3+1}=b_2,r_3\geqslant 1}}\binom{\alpha_1}{a}\binom{r_2+s_1-1}{r_2}\binom{r_4+s_2-1}{r_4}\\
&\; \times x^{\alpha_1}y\cdots x^{\alpha_{r_1}}yx^{\alpha_{r_1+1}}y^{r_2+s_1+1}x^{\widetilde{\alpha}_2}y\cdots x^{\widetilde{\alpha}_{r_3}}yx^{\widetilde{\alpha}_{r_3+1}}y^{r_4+s_2}.
\end{align*}
The second sum in the right-hand side of the above equation is
\begin{align*}
&\sum\limits_{{r_1+r_2+r_3+r_4=r,r_2\geqslant 1}\atop {\alpha_1+\cdots+\alpha_{r_1+1}=a+b_1\atop \widetilde{\alpha}_1+\cdots+\widetilde{\alpha}_{r_3+1}=b_2}}\binom{\alpha_1}{a}\binom{r_2+s_1-2}{r_2-1}\binom{r_4+s_2-1}{r_4}\\
&\; \times x^{\alpha_1}y\cdots x^{\alpha_{r_1}}yx^{\alpha_{r_1+1}}y^{r_2+s_1}x^{\widetilde{\alpha}_1}y\cdots x^{\widetilde{\alpha}_{r_3}}y x^{\widetilde{\alpha}_{r_3+1}}y^{r_4+s_2},
\end{align*}
which implies that
\begin{align*}
\Sigma_1=&\sum\limits_{{r_1+r_2+r_3+r_4=r,r_2\geqslant 1}\atop {\alpha_1+\cdots+\alpha_{r_1+1}=a+b_1\atop \widetilde{\alpha}_1+\cdots+\widetilde{\alpha}_{r_3+1}=b_2}}\binom{\alpha_1}{a}\binom{r_2+s_1-2}{r_2}\binom{r_4+s_2-1}{r_4}\\
&\; \times x^{\alpha_1}y\cdots x^{\alpha_{r_1}}yx^{\alpha_{r_1+1}}y^{r_2+s_1}x^{\widetilde{\alpha}_1}y\cdots x^{\widetilde{\alpha}_{r_3}}yx^{\widetilde{\alpha}_{r_3+1}}y^{r_4+s_2}\\
&\;+\sum\limits_{{r_1+r_3+r_4=r}\atop {\alpha_1+\cdots+\alpha_{r_1+1}=a+b_1\atop \widetilde{\alpha}_1+\cdots+\widetilde{\alpha}_{r_3+1}=b_2}}\binom{\alpha_1}{a}\binom{r_4+s_2-1}{r_4}\\
&\; \times x^{\alpha_1}y\cdots x^{\alpha_{r_1}}yx^{\alpha_{r_1+1}}y^{s_1}x^{\widetilde{\alpha}_1}y\cdots x^{\widetilde{\alpha}_{r_3}}yx^{\widetilde{\alpha}_{r_3+1}}y^{r_4+s_2}\\
=&\sum\limits_{{r_1+r_2+r_3+r_4=r}\atop {\alpha_1+\cdots+\alpha_{r_1+1}=a+b_1\atop \widetilde{\alpha}_1+\cdots+\widetilde{\alpha}_{r_3+1}=b_2}}\binom{\alpha_1}{a}\binom{r_2+s_1-2}{r_2}\binom{r_4+s_2-1}{r_4}\\
&\; \times x^{\alpha_1}y\cdots x^{\alpha_{r_1}}yx^{\alpha_{r_1+1}}y^{r_2+s_1}x^{\widetilde{\alpha}_1}y\cdots x^{\widetilde{\alpha}_{r_3}}yx^{\widetilde{\alpha}_{r_3+1}}y^{r_4+s_2}.
\end{align*}
Under the condition $r_1\geqslant 1$, we have
\begin{align}
\Sigma_1=&\sum\limits_{{r_1+r_2+r_3+r_4=r,r_1\geqslant 1}\atop {\alpha_1+\cdots+\alpha_{r_1+1}=a+b_1\atop \widetilde{\alpha}_1+\cdots+\widetilde{\alpha}_{r_3+1}=b_2}}\binom{\alpha_1}{a}\binom{r_2+s_1-2}{r_2}\binom{r_4+s_2-1}{r_4}\nonumber\\
&\; \times x^{\alpha_1}y\cdots x^{\alpha_{r_1}}yx^{\alpha_{r_1+1}}y^{r_2+s_1}x^{\widetilde{\alpha}_1}y\cdots x^{\widetilde{\alpha}_{r_3}}yx^{\widetilde{\alpha}_{r_3+1}}y^{r_4+s_2}\nonumber\\
&\;+\sum\limits_{r_1+r_2+r_3=r\atop \widetilde{\alpha}_1+\cdots+\widetilde{\alpha}_{r_2+1}=b_2}\binom{a+b_1}{a}\binom{r_1+s_1-2}{r_1}\binom{r_3+s_2-1}{r_3}\nonumber\\
&\; \times x^{a+b_1}y^{r_1+s_1}x^{\widetilde{\alpha}_1}y\cdots x^{\widetilde{\alpha}_{r_2}}yx^{\widetilde{\alpha}_{r_2+1}}y^{r_3+s_2}.
\label{Eq:Sigma-1}
\end{align}

For $\Sigma_2$, we have
\begin{align*}
\Sigma_2=&\sum\limits_{{{1\leqslant l\leqslant s_1\atop 0\leqslant a_1\leqslant a-1}\atop r_1+r_2+r_3=r}\atop{\alpha_1+\cdots+\alpha_{l+1}=a+b_1,\alpha_1\geqslant a_1+b_1,\alpha_{l+1}\geqslant 1\atop \widetilde{\alpha}_1+\cdots+\widetilde{\alpha}_{r_2+1}=b_2,\widetilde{\alpha}_1\geqslant 1}}\binom{a_1+b_1-1}{b_1-1}\binom{r_1+s_1-l}{r_1}\binom{r_3+s_2-1}{r_3}\\
&\;\times x^{\alpha_1}y\cdots x^{\alpha_l}yx^{\alpha_{l+1}}y^{r_1+s_1-l}x^{\widetilde{\alpha}_1}y\cdots x^{\widetilde{\alpha}_{r_2}}yx^{\widetilde{\alpha}_{r_2+1}}y^{r_3+s_2}\\
=&\sum\limits_{{1\leqslant l\leqslant s_1\atop r_1+r_2+r_3=r}\atop{\alpha_1+\cdots+\alpha_{l+1}=a+b_1,\alpha_{l+1}\geqslant 1\atop \widetilde{\alpha}_1+\cdots+\widetilde{\alpha}_{r_2+1}=b_2,\widetilde{\alpha}_1\geqslant 1}}\binom{\alpha_1}{b_1}\binom{r_1+s_1-l}{r_1}\binom{r_3+s_2-1}{r_3}\\
&\;\times x^{\alpha_1}y\cdots x^{\alpha_l}yx^{\alpha_{l+1}}y^{r_1+s_1-l}x^{\widetilde{\alpha}_1}y\cdots x^{\widetilde{\alpha}_{r_2}}yx^{\widetilde{\alpha}_{r_2+1}}y^{r_3+s_2}.
\end{align*}
Neglecting the conditions $\widetilde{\alpha}_1\geqslant 1$, we get
\begin{align*}
\Sigma_2=&\sum\limits_{{1\leqslant l\leqslant s_1\atop r_1+r_2+r_3=r}\atop{\alpha_1+\cdots+\alpha_{l+1}=a+b_1,\alpha_{l+1}\geqslant 1\atop \widetilde{\alpha}_1+\cdots+\widetilde{\alpha}_{r_2+1}=b_2}}\binom{\alpha_1}{b_1}\binom{r_1+s_1-l}{r_1}\binom{r_3+s_2-1}{r_3}\\
&\;\times x^{\alpha_1}y\cdots x^{\alpha_l}yx^{\alpha_{l+1}}y^{r_1+s_1-l}x^{\widetilde{\alpha}_1}y\cdots x^{\widetilde{\alpha}_{r_2}}yx^{\widetilde{\alpha}_{r_2+1}}y^{r_3+s_2}\\
&\;-\sum\limits_{{1\leqslant l\leqslant s_1\atop r_1+r_2+r_3=r,r_1\geqslant 1}\atop{\alpha_1+\cdots+\alpha_{l+1}=a+b_1,\alpha_{l+1}\geqslant 1\atop \widetilde{\alpha}_1+\cdots+\widetilde{\alpha}_{r_2+1}=b_2}}\binom{\alpha_1}{b_1}\binom{r_1+s_1-l-1}{r_1-1}\binom{r_3+s_2-1}{r_3}\\
&\;\times x^{\alpha_1}y\cdots x^{\alpha_l}yx^{\alpha_{l+1}}y^{r_1+s_1-l}x^{\widetilde{\alpha}_1}y\cdots x^{\widetilde{\alpha}_{r_2}}yx^{\widetilde{\alpha}_{r_2+1}}y^{r_3+s_2}\\
=&\sum\limits_{{1\leqslant l\leqslant s_1\atop r_1+r_2+r_3=r,r_1\geqslant 1}\atop{\alpha_1+\cdots+\alpha_{l+1}=a+b_1,\alpha_{l+1}\geqslant 1\atop \widetilde{\alpha}_1+\cdots+\widetilde{\alpha}_{r_2+1}=b_2}}\binom{\alpha_1}{b_1}\binom{r_1+s_1-l-1}{r_1}\binom{r_3+s_2-1}{r_3}\\
&\;\times x^{\alpha_1}y\cdots x^{\alpha_l}yx^{\alpha_{l+1}}y^{r_1+s_1-l}x^{\widetilde{\alpha}_1}y\cdots x^{\widetilde{\alpha}_{r_2}}yx^{\widetilde{\alpha}_{r_2+1}}y^{r_3+s_2}\\
&\;+\sum\limits_{{1\leqslant l\leqslant s_1\atop r_1+r_2=r}\atop{\alpha_1+\cdots+\alpha_{l+1}=a+b_1,\alpha_{l+1}\geqslant 1\atop \widetilde{\alpha}_1+\cdots+\widetilde{\alpha}_{r_1+1}=b_2}}\binom{\alpha_1}{b_1}\binom{r_2+s_2-1}{r_2}\\
&\;\times x^{\alpha_1}y\cdots x^{\alpha_l}yx^{\alpha_{l+1}}y^{s_1-l}x^{\widetilde{\alpha}_1}y\cdots x^{\widetilde{\alpha}_{r_1}}yx^{\widetilde{\alpha}_{r_1+1}}y^{r_2+s_2}.
\end{align*}
We denote the first sum in the right-hand side of the above equation by $\Sigma_{21}$, and the second sum by $\Sigma_{22}$. Then without the condition $\alpha_{l+1}\geqslant 1$, we have
\begin{align*}
\Sigma_{21}=&\sum\limits_{{1\leqslant l\leqslant s_1\atop r_1+r_2+r_3=r,r_1\geqslant 1}\atop{\alpha_1+\cdots+\alpha_{l+1}=a+b_1\atop \widetilde{\alpha}_1+\cdots+\widetilde{\alpha}_{r_2+1}=b_2}}\binom{\alpha_1}{b_1}\binom{r_1+s_1-l-1}{r_1}\binom{r_3+s_2-1}{r_3}\\
&\;\times x^{\alpha_1}y\cdots x^{\alpha_l}yx^{\alpha_{l+1}}y^{r_1+s_1-l}x^{\widetilde{\alpha}_1}y\cdots x^{\widetilde{\alpha}_{r_2}}yx^{\widetilde{\alpha}_{r_2+1}}y^{r_3+s_2}\\
&\;-\sum\limits_{{0\leqslant l\leqslant s_1-1\atop r_1+r_2+r_3=r,r_1\geqslant 1}\atop{\alpha_1+\cdots+\alpha_{l+1}=a+b_1\atop \widetilde{\alpha}_1+\cdots+\widetilde{\alpha}_{r_2+1}=b_2}}\binom{\alpha_1}{b_1}\binom{r_1+s_1-l-2}{r_1}\binom{r_3+s_2-1}{r_3}\\
&\;\times x^{\alpha_1}y\cdots x^{\alpha_l}yx^{\alpha_{l+1}}y^{r_1+s_1-l}x^{\widetilde{\alpha}_1}y\cdots x^{\widetilde{\alpha}_{r_2}}yx^{\widetilde{\alpha}_{r_2+1}}y^{r_3+s_2}\\
=&\sum\limits_{{1\leqslant l\leqslant s_1-1\atop r_1+r_2+r_3=r,r_1\geqslant 1}\atop{\alpha_1+\cdots+\alpha_{l+1}=a+b_1\atop \widetilde{\alpha}_1+\cdots+\widetilde{\alpha}_{r_2+1}=b_2}}\binom{\alpha_1}{b_1}\binom{r_1+s_1-l-2}{r_1-1}\binom{r_3+s_2-1}{r_3}\\
&\;\times x^{\alpha_1}y\cdots x^{\alpha_l}yx^{\alpha_{l+1}}y^{r_1+s_1-l}x^{\widetilde{\alpha}_1}y\cdots x^{\widetilde{\alpha}_{r_2}}yx^{\widetilde{\alpha}_{r_2+1}}y^{r_3+s_2}\\
&\;-\sum\limits_{{r_1+r_2+r_3=r,r_1\geqslant 1}\atop {\widetilde{\alpha}_1+\cdots+\widetilde{\alpha}_{r_2+1}=b_2}}\binom{a+b_1}{b_1}\binom{r_1+s_1-2}{r_1}\binom{r_3+s_2-1}{r_3}\\
&\;\times x^{a+b_1}y^{r_1+s_1}x^{\widetilde{\alpha}_1}y\cdots x^{\widetilde{\alpha}_{r_2}}yx^{\widetilde{\alpha}_{r_2+1}}y^{r_3+s_2},
\end{align*}
and
\begin{align*}
\Sigma_{22}=&\sum\limits_{{1\leqslant l\leqslant s_1\atop r_1+r_2=r}\atop{\alpha_1+\cdots+\alpha_{l+1}=a+b_1\atop \widetilde{\alpha}_1+\cdots+\widetilde{\alpha}_{r_1+1}=b_2}}\binom{\alpha_1}{b_1}\binom{r_2+s_2-1}{r_2}\\
&\;\times x^{\alpha_1}y\cdots x^{\alpha_l}yx^{\alpha_{l+1}}y^{s_1-l}x^{\widetilde{\alpha}_1}y\cdots x^{\widetilde{\alpha}_{r_1}}yx^{\widetilde{\alpha}_{r_1+1}}y^{r_2+s_2}\\
&\;-\sum\limits_{{0\leqslant l\leqslant s_1-1\atop r_1+r_2=r}\atop{\alpha_1+\cdots+\alpha_{l+1}=a+b_1\atop \widetilde{\alpha}_1+\cdots+\widetilde{\alpha}_{r_1+1}=b_2}}\binom{\alpha_1}{b_1}\binom{r_2+s_2-1}{r_2}\\
&\;\times x^{\alpha_1}y\cdots x^{\alpha_l}yx^{\alpha_{l+1}}y^{s_1-l}x^{\widetilde{\alpha}_1}y\cdots x^{\widetilde{\alpha}_{r_1}}yx^{\widetilde{\alpha}_{r_1+1}}y^{r_2+s_2}\\
=&\sum\limits_{{r_1+r_2=r}\atop{\alpha_1+\cdots+\alpha_{s_1+1}=a+b_1\atop \widetilde{\alpha}_1+\cdots+\widetilde{\alpha}_{r_1+1}=b_2}}\binom{\alpha_1}{b_1}\binom{r_2+s_2-1}{r_2}\\
&\;\times x^{\alpha_1}y\cdots x^{\alpha_{s_1}}yx^{\alpha_{s_1+1}+\widetilde{\alpha}_1}yx^{\widetilde{\alpha}_2}y\cdots x^{\widetilde{\alpha}_{r_1}}yx^{\widetilde{\alpha}_{r_1+1}}y^{r_2+s_2}\\
&\;-\sum\limits_{{r_1+r_2=r}\atop{\widetilde{\alpha}_1+\cdots+\widetilde{\alpha}_{r_1+1}=b_2}}\binom{a_1+b_1}{b_1}\binom{r_2+s_2-1}{r_2}\\
&\;\times x^{a+b_1}y^{s_1}x^{\widetilde{\alpha}_1}y\cdots x^{\widetilde{\alpha}_{r_1}}yx^{\widetilde{\alpha}_{r_1+1}}y^{r_2+s_2}.
\end{align*}
Hence we get
\begin{align}
\Sigma_2=&\sum\limits_{{1\leqslant l\leqslant s_1-1\atop r_1+r_2+r_3=r,r_1\geqslant 1}\atop{\alpha_1+\cdots+\alpha_{l+1}=a+b_1\atop \widetilde{\alpha}_1+\cdots+\widetilde{\alpha}_{r_2+1}=b_2}}\binom{\alpha_1}{b_1}\binom{r_1+s_1-l-2}{r_1-1}\binom{r_3+s_2-1}{r_3}\nonumber\\
&\;\times x^{\alpha_1}y\cdots x^{\alpha_l}yx^{\alpha_{l+1}}y^{r_1+s_1-l}x^{\widetilde{\alpha}_1}y\cdots x^{\widetilde{\alpha}_{r_2}}yx^{\widetilde{\alpha}_{r_2+1}}y^{r_3+s_2}\nonumber\\
&\;-\sum\limits_{{r_1+r_2+r_3=r}\atop {\widetilde{\alpha}_1+\cdots+\widetilde{\alpha}_{r_2+1}=b_2}}\binom{a+b_1}{b_1}\binom{r_1+s_1-2}{r_1}\binom{r_3+s_2-1}{r_3}\nonumber\\
&\;\times x^{a+b_1}y^{r_1+s_1}x^{\widetilde{\alpha}_1}y\cdots x^{\widetilde{\alpha}_{r_2}}yx^{\widetilde{\alpha}_{r_2+1}}y^{r_3+s_2}\nonumber\\
&\;+\sum\limits_{{r_1+r_2=r}\atop{\alpha_1+\cdots+\alpha_{s_1}+\beta+\widetilde{\alpha}_2+\cdots+\widetilde{\alpha}_{r_1+1}=a+b_1+b_2\atop a+b_1\geqslant \alpha_1+\cdots+\alpha_{s_1}\geqslant a+b_1-\beta}}\binom{\alpha_1}{b_1}\binom{r_2+s_2-1}{r_2}\nonumber\\
&\;\times x^{\alpha_1}y\cdots x^{\alpha_{s_1}}yx^{\beta}yx^{\widetilde{\alpha}_2}y\cdots x^{\widetilde{\alpha}_{r_1}}yx^{\widetilde{\alpha}_{r_1+1}}y^{r_2+s_2}.
\label{Eq:Sigma-2}
\end{align}

Similarly, for $\Sigma_3$, we have
\begin{align*}
\Sigma_3=&\sum\limits_{{{1\leqslant k\leqslant b_2\atop a_1+a_3\leqslant a-1}\atop{r_1+r_2=r\atop \alpha_1+\cdots+\alpha_{s_1+1}=a+b_1-a_3}}\atop {\widetilde{\alpha}_1+\cdots+\widetilde{\alpha}_{r_1+1}=a_3+b_2\atop \alpha_1\geqslant a_1+b_1,\alpha_{s_1+1}\geqslant 1,\widetilde{\alpha}_1\geqslant a_3+k}}\binom{a_1+b_1-1}{b_1-1}\binom{a_3+k-1}{k-1}\binom{r_2+s_2-1}{r_2}\\
&\;\times x^{\alpha_1}y\cdots x^{\alpha_{s_1}}yx^{\alpha_{s_1+1}+\widetilde{\alpha}_1}yx^{\widetilde{\alpha}_2}y\cdots x^{\widetilde{\alpha}_{r_1}}yx^{\widetilde{\alpha}_{r_1+1}}y^{r_2+s_2}\\
=&\sum\limits_{{{a_1+a_3\leqslant a-1}\atop{r_1+r_2=r\atop \alpha_1+\cdots+\alpha_{s_1+1}=a+b_1-a_3}}\atop {\widetilde{\alpha}_1+\cdots+\widetilde{\alpha}_{r_1+1}=a_3+b_2\atop \alpha_1\geqslant a_1+b_1,\alpha_{s_1+1}\geqslant 1,\widetilde{\alpha}_1\geqslant a_3+1}}\binom{a_1+b_1-1}{b_1-1}\binom{\widetilde{\alpha}_1}{a_3+1}\binom{r_2+s_2-1}{r_2}\\
&\;\times x^{\alpha_1}y\cdots x^{\alpha_{s_1}}yx^{\alpha_{s_1+1}+\widetilde{\alpha}_1}yx^{\widetilde{\alpha}_2}y\cdots x^{\widetilde{\alpha}_{r_1}}yx^{\widetilde{\alpha}_{r_1+1}}y^{r_2+s_2}\\
=&\sum\limits_{{{r_1+r_2=r\atop \alpha_1+\cdots+\alpha_{s_1+1}=a+b_1-a_3}}\atop {\widetilde{\alpha}_1+\cdots+\widetilde{\alpha}_{r_1+1}=a_3+b_2\atop \alpha_{s_1+1}\geqslant 1,\widetilde{\alpha}_1\geqslant a_3+1}}\binom{\alpha_1}{b_1}\binom{\widetilde{\alpha}_1}{a_3+1}\binom{r_2+s_2-1}{r_2}\\
&\;\times x^{\alpha_1}y\cdots x^{\alpha_{s_1}}yx^{\alpha_{s_1+1}+\widetilde{\alpha}_1}yx^{\widetilde{\alpha}_2}y\cdots x^{\widetilde{\alpha}_{r_1}}yx^{\widetilde{\alpha}_{r_1+1}}y^{r_2+s_2}\\
=&\sum\limits_{{r_1+r_2=r\atop a_3\leqslant a+b_1-\alpha_1-\cdots-\alpha_{s_1}-1}\atop {\alpha_1+\cdots+\alpha_{s_1}\geqslant a+b_1-\beta+1\atop\alpha_1+\cdots+\alpha_{s_1}+\beta+\widetilde{\alpha}_2+\cdots+\widetilde{\alpha}_{r_1+1}=a+b_1+b_2}}\binom{\alpha_1}{b_1}
\binom{\beta+\alpha_1+\cdots+\alpha_{s_1}-a-b_1+a_3}{a_3+1}\\
&\;\times\binom{r_2+s_2-1}{r_2}x^{\alpha_1}y\cdots x^{\alpha_{s_1}}yx^{\beta}yx^{\widetilde{\alpha}_2}y\cdots x^{\widetilde{\alpha}_{r_1}}yx^{\widetilde{\alpha}_{r_1+1}}y^{r_2+s_2}\\
=&\sum\limits_{{r_1+r_2=r\atop a+b_1\geqslant\alpha_1+\cdots+\alpha_{s_1}\geqslant a+b_1-\beta}\atop  \alpha_1+\cdots+\alpha_{s_1}+\beta+\widetilde{\alpha}_2+\cdots+\widetilde{\alpha}_{r_1+1}=a+b_1+b_2}\binom{\alpha_1}{b_1}
\left[\binom{\beta}{a+b_1-\alpha_1-\cdots-\alpha_{s_1}}-1\right]\\
&\;\times\binom{r_2+s_2-1}{r_2}x^{\alpha_1}y\cdots x^{\alpha_{s_1}}yx^{\beta}yx^{\widetilde{\alpha}_2}y\cdots x^{\widetilde{\alpha}_{r_1}}yx^{\widetilde{\alpha}_{r_1+1}}y^{r_2+s_2}\\
=&\sum\limits_{{r_1+r_2=r}\atop  \alpha_1+\cdots+\alpha_{s_1}+\beta+\widetilde{\alpha}_2+\cdots+\widetilde{\alpha}_{r_1+1}=a+b_1+b_2}\binom{\alpha_1}{b_1}
\binom{\beta}{a+b_1-\alpha_1-\cdots-\alpha_{s_1}}\\
&\;\times\binom{r_2+s_2-1}{r_2}x^{\alpha_1}y\cdots x^{\alpha_{s_1}}yx^{\beta}yx^{\widetilde{\alpha}_2}y\cdots x^{\widetilde{\alpha}_{r_1}}yx^{\widetilde{\alpha}_{r_1+1}}y^{r_2+s_2}\\
&\;-\sum\limits_{{r_1+r_2=r\atop a+b_1\geqslant \alpha_1+\cdots+\alpha_{s_1}\geqslant a+b_1-\beta}\atop  \alpha_1+\cdots+\alpha_{s_1}+\beta+\widetilde{\alpha}_2+\cdots+\widetilde{\alpha}_{r_1+1}=a+b_1+b_2}\binom{\alpha_1}{b_1}\binom{r_2+s_2-1}{r_2}\\
&\;\times x^{\alpha_1}y\cdots x^{\alpha_{s_1}}yx^{\beta}yx^{\widetilde{\alpha}_2}y\cdots x^{\widetilde{\alpha}_{r_1}}yx^{\widetilde{\alpha}_{r_1+1}}y^{r_2+s_2}.
\end{align*}
Under the condition $r_1\geqslant 1$ in the first sum in the above equation, we have
\begin{align}
\Sigma_3=&\sum\limits_{{r_1+r_2=r,r_1\geqslant 1}\atop  \alpha_1+\cdots+\alpha_{s_1}+\beta+\widetilde{\alpha}_2+\cdots+\widetilde{\alpha}_{r_1+1}=a+b_1+b_2}\binom{\alpha_1}{b_1}
\binom{\beta}{a+b_1-\alpha_1-\cdots-\alpha_{s_1}}\nonumber\\
&\;\times\binom{r_2+s_2-1}{r_2}x^{\alpha_1}y\cdots x^{\alpha_{s_1}}yx^{\beta}yx^{\widetilde{\alpha}_2}y\cdots x^{\widetilde{\alpha}_{r_1}}yx^{\widetilde{\alpha}_{r_1+1}}y^{r_2+s_2}\nonumber\\
&\;+\sum\limits_{ \alpha_1+\cdots+\alpha_{s_1}+\beta=a+b_1+b_2}\binom{\alpha_1}{b_1}
\binom{\beta}{a+b_1-\alpha_1-\cdots-\alpha_{s_1}}\nonumber\\
&\;\times\binom{r+s_2-1}{r}x^{\alpha_1}y\cdots x^{\alpha_{s_1}}yx^{\beta}y^{r+s_2}\nonumber\\
&\;-\sum\limits_{{r_1+r_2=r\atop a+b_1\geqslant\alpha_1+\cdots+\alpha_{s_1}\geqslant a+b_1-\beta}\atop  \alpha_1+\cdots+\alpha_{s_1}+\beta+\widetilde{\alpha}_2+\cdots+\widetilde{\alpha}_{r_1+1}=a+b_1+b_2}\binom{\alpha_1}{b_1}\binom{r_2+s_2-1}{r_2}\nonumber\\
&\;\times x^{\alpha_1}y\cdots x^{\alpha_{s_1}}yx^{\beta}yx^{\widetilde{\alpha}_2}y\cdots x^{\widetilde{\alpha}_{r_1}}yx^{\widetilde{\alpha}_{r_1+1}}y^{r_2+s_2}.
\label{Eq:Sigma-3}
\end{align}

Finally, for $\Sigma_4$, we have
\begin{align*}
\Sigma_4=&\sum\limits_{{{1\leqslant l\leqslant s_2\atop a_1+a_2+a_3+a_4=a-1}\atop{\alpha_1+\cdots+\alpha_{s_1}=a_1+a_2+b_1\atop \widetilde{\alpha}_1+\cdots+\widetilde{\alpha}_{l+1}=a_3+a_4+b_2+1}}\atop \alpha_1\geqslant a_1+b_1,\widetilde{\alpha}_1\geqslant a_3+b_2,\widetilde{\alpha}_{l+1}\geqslant 1}\binom{a_1+b_1-1}{b_1-1}\binom{a_3+b_2-1}{b_2-1}\binom{r+s_2-l}{r}\\
&\;\times x^{\alpha_1}y\cdots x^{\alpha_{s_1}}yx^{\widetilde{\alpha}_1}y\cdots x^{\widetilde{\alpha}_l}yx^{\widetilde{\alpha}_{l+1}}y^{r+s_2-l}\\
=&\sum\limits_{{{1\leqslant l\leqslant s_2\atop \alpha_1+\cdots+\alpha_{s_1}+\widetilde{\alpha}_1+\cdots+\widetilde{\alpha}_{l+1}=a+b_1+b_2}\atop{a_1\leqslant \alpha_1+\cdots+\alpha_{s_1}-b_1\atop a_3\leqslant \widetilde{\alpha}_1+\cdots+\widetilde{\alpha}_{l+1}-b_2-1}}\atop \alpha_1\geqslant a_1+b_1,\widetilde{\alpha}_1\geqslant a_3+b_2,\widetilde{\alpha}_{l+1}\geqslant 1}\binom{a_1+b_1-1}{b_1-1}\binom{a_3+b_2-1}{b_2-1}\binom{r+s_2-l}{r}\\
&\;\times x^{\alpha_1}y\cdots x^{\alpha_{s_1}}yx^{\widetilde{\alpha}_1}y\cdots x^{\widetilde{\alpha}_l}yx^{\widetilde{\alpha}_{l+1}}y^{r+s_2-l}\\
=&\sum\limits_{{\alpha_1+\cdots+\alpha_{s_1}+\widetilde{\alpha}_1+\cdots+\widetilde{\alpha}_{l+1}=a+b_1+b_2}\atop 1\leqslant l\leqslant s_2,\widetilde{\alpha}_{l+1}\geqslant 1}\binom{\alpha_1}{b_1}\binom{\widetilde{\alpha}_1}{b_2}\binom{r+s_2-l}{r}\\
&\;\times x^{\alpha_1}y\cdots x^{\alpha_{s_1}}yx^{\widetilde{\alpha}_1}y\cdots x^{\widetilde{\alpha}_l}yx^{\widetilde{\alpha}_{l+1}}y^{r+s_2-l}.
\end{align*}
Without the condition $\widetilde{\alpha}_{l+1}\geqslant 1$, we get
\begin{align}
\Sigma_4=&\sum\limits_{{\alpha_1+\cdots+\alpha_{s_1}+\widetilde{\alpha}_1+\cdots+\widetilde{\alpha}_{l+1}=a+b_1+b_2}\atop 1\leqslant l\leqslant s_2}\binom{\alpha_1}{b_1}\binom{\widetilde{\alpha}_1}{b_2}\binom{r+s_2-l-1}{r-1}\nonumber\\
&\;\times x^{\alpha_1}y\cdots x^{\alpha_{s_1}}yx^{\widetilde{\alpha}_1}y\cdots x^{\widetilde{\alpha}_l}yx^{\widetilde{\alpha}_{l+1}}y^{r+s_2-l}\nonumber\\
&-\sum\limits_{\alpha_1+\cdots+\alpha_{s_1}+\widetilde{\alpha}_1=a+b_1+b_2}\binom{\alpha_1}{b_1}\binom{\widetilde{\alpha}_1}{b_2}\binom{r+s_2-1}{r}
\nonumber\\
&\;\times x^{\alpha_1}y\cdots x^{\alpha_{s_1}}yx^{\widetilde{\alpha}_1}y^{r+s_2}.
\label{Eq:Sigma-4}
\end{align}

The equations \eqref{Eq:Sigma-1}-\eqref{Eq:Sigma-4} imply that the sum $\Sigma_1+\Sigma_2+\Sigma_3+\Sigma_4$ is just the right-hand side of the formula \eqref{Eq:Shuffle-Res-1-2}. Then we proved that the formula \eqref{Eq:Shuffle-Res-1-2-Old} can be deduced from \eqref{Eq:Shuffle-Res-1-2}.


\end{document}